\documentclass[11pt]{article}
\usepackage{graphicx}
\usepackage{epsfig}
\usepackage{latexsym}
\newsavebox{\toy}
\savebox{\toy}{\framebox[0.65em]{\rule{0cm}{1ex}}}
\newcommand{\QED}{\usebox{\toy}\end{demo}}
\newenvironment{property}%
{\begin{list}{}{\setlength{\rightmargin}{0pt}%
\setlength{\itemsep}{0pt}}}{\end{list}}
\newlength{\templength}
\newcommand{\bp}{\setlength{\templength}{\labelwidth}%
\setlength{\labelwidth}{2em}\begin{property}}
\newcommand{\ep}{\end{property}\setlength{\labelwidth}{\templength}}
\newtheorem{theorem}{Theorem}[subsection]
\newtheorem{lemma}[theorem]{Lemma}
\newtheorem{proposition}[theorem]{Proposition}
\newtheorem{corollary}[theorem]{Corollary}
\newtheorem{assumption}{Assumption}
\newtheorem{definition}{Definition}[subsection]
\newtheorem{remark}{Remark}[subsection]
\newtheorem{conjecture}{Conjecture}[subsection]
\newtheorem{exercise}{Exercise}[subsection]

\newcommand{\Thm}[1]{Theorem \ref{Thm.#1}}
\newcommand{\Lem}[1]{Lemma \ref{Lem.#1}}

\newcommand{\Prop}[1]{Proposition \ref{Prop.#1}}

\newcommand{\Theorem}[1]{\begin{theorem}\label{Thm.#1}}
\newcommand{\Lemma}[1]{\begin{lemma}\label{Lem.#1}}
\newcommand{\Proposition}[1]{\begin{proposition}\label{Prop.#1}}
\newcommand{\Corollary}[1]{\begin{corollary}\label{Cor.#1}}
\newcommand{\Assumption}[1]{\begin{assumption}\label{Ass.#1}\rm}
\newcommand{\Definition}[1]{\begin{definition}\label{Def.#1}\rm}
\newcommand{\Remark}[1]{\begin{remark}\label{Rem.#1}\rm }
\newcommand{\Exercise}[1]{\begin{exercise}\label{Exe.#1}\rm }

\newcommand{\bd}{\begin{displaymath}}
\newcommand{\ed}{\end{displaymath}}
\newcommand{\bdn}{\begin{equation}}
\newcommand{\bdnl}{\begin{equation}\label}
\newcommand{\edn}{\end{equation}}
\newcommand{\barray}{\begin{array}}
\newcommand{\earray}{\end{array}}
\newcommand{\bds}{\begin{description}}
\newcommand{\eds}{\end{description}}
\newcommand{\bitemize}{\begin{itemize}}
\newcommand{\eitemize}{\end{itemize}}
\newcommand{\benumerate}{\begin{enumerate}}
\newcommand{\eenumerate}{\end{enumerate}}
\newcommand{\btabbing}{\begin{tabbing}}
\newcommand{\etabbing}{\end{tabbing}}
\newcommand{\bcenter}{\begin{center}}
\newcommand{\ecenter}{\end{center}}
\newcommand{\bflushright}{\begin{flushright}}
\newcommand{\bflushleft}{\begin{flushleft}}
\newcommand{\eflushright}{\end{flushright}}
\newcommand{\eflushleft}{\end{flushleft}}
\newcommand{\bdnn }{\begin{eqnarray*}}
\newcommand{\ednn }{\end{eqnarray*}}
\newcommand{\bdmn}{\begin{eqnarray}}
\newcommand{\edmn}{\end{eqnarray}}

\newcommand{\nn}{\nonumber}
\newcommand{\SSC}[1]{\section{#1}\setcounter{equation}{0}}

\newcounter{biblio}
\newenvironment{references}%
{\begin{list}{[\arabic{biblio}]}{\usecounter{biblio}%
\setlength{\leftmargin}{2.5em}\setlength{\rightmargin}{0pt}%
\setlength{\labelwidth}{2em}\setlength{\itemsep}{0pt}}}{\end{list}}
\newcommand{\References}%
{\vspace{2.8ex plus .3ex minus .3ex}%
\begin{center}{\bf References}\end{center}\begin{references}}

\usepackage{amssymb}
\usepackage{mathrsfs}
\newcommand{\N}{{\mathbb{N}}}
\newcommand{\Z}{{\mathbb{Z}}}
\newcommand{\zd}{\Z^d}
\newcommand{\Q}{{\mathbb{Q}}}

\newcommand{\R}{{\mathbb{R}}}
\newcommand{\rd}{\R^d}

\newcommand{\ra }{\rightarrow }

\newcommand{\vvs}{\vspace{2ex}}
\newcommand{\vs}{\vspace{1ex}}

\newcommand{\lef}{\left}
\newcommand{\rig}{\right}
\newcommand{\ri}{\right}
\newcommand{\st}{\stackrel}

\newcommand{\8}{\infty}

\newcommand{\dps}{\displaystyle}
\newcommand{\sub}{\subset}

\newcommand{\bsh}{\backslash}
\newcommand{\pri}{\prime}

\newcommand{\inflim}{\mathop{\underline{\lim}}}

\newcommand{\epty}{\emptyset}
\renewcommand{\a}{\alpha}

\newcommand{\del}{\delta}
\newcommand{\D}{\Delta}

\newcommand{\e}{\varepsilon}

\newcommand{\z}{\zeta}
\newcommand{\h}{\eta}
\newcommand{\tht}{\theta}

\newcommand{\lm}{\lambda}

\newcommand{\m}{\mu}
\newcommand{\n}{\nu}

\newcommand{\s}{\sigma}

\newcommand{\w}{\omega}

\newcommand{\W}{\Omega}

\newcommand{\cF }{{\cal F}}
\newcommand{\cG }{{\cal G}}

\newcommand{\cP }{{\cal P}}

\newcommand{\cT }{{\cal T}}

\makeatletter
\def\section{\@startsection{section}{1}{\z@}{-3.5ex plus -1ex minus
 -.2ex}{2.3ex plus .2ex}{\bf}}
\def\subsection{\@startsection{subsection}{2}{\z@}{-3.25ex plus -1ex minus
 -.2ex}{1.5ex plus .2ex}{\bf}}
\makeatother

\begin{document}
\parindent=0pt

\bcenter

\large{\bf Branching Random Walks in Space-Time Random Environment: 
Survival Probability, Global and Local Growth Rates}\footnote{\today}

\vvs \normalsize

\noindent Francis Comets\footnote{Partially supported by ANR Polintbio} 
and Nobuo Yoshida\footnote{Partially supported by JSPS 
Grant-in-Aid for Scientific Research, Kiban (C) 21540125} \\

\ecenter

\begin{abstract}
We study the survival probability and the 
growth rate for  branching random walks in random environment (BRWRE). 
The particles perform simple symmetric random 
walks on the $d$-dimensional integer lattice, while at each 
time unit, they split into independent copies according 
to time-space i.i.d. offspring distributions. 
The BRWRE is naturally associated with the 
directed polymers in random environment (DPRE), 
for which the quantity called the free energy is well studied. 
We discuss the survival 
probability (both global and local) 
for BRWRE and give a criterion for its positivity
in terms of the free energy of the associated DPRE.
We also show that the global growth rate for the number of particles 
in BRWRE is given by  the free energy of 
the associated DPRE, though the local growth rate
is  given by  the directional free energy.
\end{abstract}

\vspace{1cm}
\footnotesize 
\noindent{\bf Short Title.} Branching Walks in Random Environment
 


\normalsize

\small

\tableofcontents

\normalsize
\SSC{Introduction}
\hspace{3mm} We consider particles in $\zd$, performing random walks 
and branching into independent copies at each step.
There is initially one particle located in the origin.
When a particle occupies a site $x \in \zd$ at time 
$t  \in \N=\{0,1,\ldots\}$, then, 
it moves to a randomly choosen  adjacent site $y$ 
at time $t+1$ and is replaced by 
$k$ new particles with probability $q_{t,x}(k)$ ($k  \in \N$). 
We assume that the offspring distributions  
$q_{t,x}=(q_{t,x}(k))_{k \in \N}$ form an i.i.d. sequence indexed by 
time $t$ and 
space $x$ variables, and we refer to this sequence as the environment.
The above model, that we call  branching random walks in (space-time) random 
environment and abreviate as BRWRE, was 
introduced by Birkner in his thesis \cite{Bir03} under the supervision of
Wakolbinger, to analyse
interactions between particles sharing the same environment,
but independent otherwise.  
Birkner, Geiger and Kersting \cite{BGK05}
obtain sufficient conditions for survival and for extinction, Hu and Yoshida
investigate the localization properties of the empirical particle density
in \cite{HuYo07}, while the second author studies in \cite{Yo08a}
the diffusive behaviour for $d \geq 3$ and when the diffusion dominates the 
disorder of the environment. See also an improvement of \cite{Yo08a} by 
Nakashima \cite{Na09}.

\hspace{3mm}
Let $B_t$ be the set of particles in BRWRE at time $t$, and denote 
the survival event by
\bdnl{DefSurv}
\{ {\rm survival}\}\st{\rm def.}{=}
\{ B_t \neq \epty \; \; \mbox{for all $t \ge 0$}\}.
\edn
The first fundamental question we address in this paper, is:
\bdnl{GWsurv}
\mbox{When is 
the probability of the above event (survival probability) positive? }
\edn
When 
it is positive indeed, 
a further natural question will be:
\bdnl{GWgrow}
\mbox{How fast does the total population $|B_t|$ grow as $t \ra \8$?}
\edn
As will be explained in section \ref{sec:dpre} below, 
the BRWRE is naturally associated with 
a model of directed polymers in random environment (DPRE) --
see e.g. 
\cite{CaHu03,CSY03,CSY04} --, which describes the population mean in a fixed
environment. We will show that questions 
(\ref{GWsurv})--(\ref{GWgrow}) intrinsically 
relate to the behavior of the associated DPRE,
and we will answer them in terms of the free energy for the 
DPRE. We will prove that survival occurs with positive probability
when the free energy is positive, while extinction occurs almost surely
when it is negative; we leave open, the case when the free energy is zero.
On the event of survival, the growth rate for the
population is the same as for its expectation given the environment.
The latter one can be strictly smaller than the one for the unconditional
expectation. Besides the overall population, we will consider also 
the set of particles moving in a specific direction. An interesting 
point is that martingale theory does not work well in this problem
for a large range of parameters (offspring distribution, environmental law,
space 
dimension, \ldots),
neither under the annealed law nor the quenched one, and we will have to 
resort to different techniques.

Other models of branching random walks in random environment are 
often considered in the literature. A first one, introduced in \cite{Sm68}, 
is when the offspring 
distribution is a random i.i.d. sequence depending on time only; we will 
call it the Smith-Wilkinson model, and it  plays a crucial role 
for comparison with ours. Another popular model is when the offspring 
distribution is a random i.i.d. sequence depending on space only;
our model relates to that one in space dimension $d+1$ by adding an 
(infinitely) large drift in the first direction,
or, equivalently, by considering crossings in that direction.
In that model, the interest for global/local growth rates started in
\cite{GdH91} in one spatial dimension, allowing explicit results, and 
in \cite{CP07} for higher dimension but in the absence of extinction;
the reader is refered to \cite{GMPV09} and \cite{Mu08} for recent 
contributions on global/local survival and recurrence/transience issues.  
We emphasize that the possibility of extinction makes the study
of growth rates much more delicate, e.g., many constructions in  \cite{CP07}
require there is at least one particle in the population at all times. 

{\bf Notations}
We write 
$\N^*=\{1,2,\ldots\}=\N \setminus \{0\}$.
For $x \in \rd$, $|x|$ denotes the $\ell^1$ norm:  
$|x|=|x_1|+\ldots +|x_d|$. 
We write $P[X]=\int X \; dP$ and 
$P[X:A]=\int_A X \; dP$ for a random variable $X$ and an 
event $A$ on some  probability space.
For two events $A,B$, we write $A \sub B$ a.s. iff
$P(A \bsh B)=0$. Similarly, $A = B$ a.s. means that 
$P(A \bsh B)=P(B \bsh A)=0$.
We denote by 
$\cP (\N)=\Big\{ {\bf q}=(q (k))_{k \in \N} \in [0,1]^\N \; : \; \sum_{k \in
\N}q (k)=1\Big\}$
the set of probability measures on $\N$.

\subsection{Backgrounds}
We begin by recalling some classical results for 
the Galton-Watson process and the Smith-Wilkinson process. 
To make the notation consistent with 
later ones, we formulate these processes as the time evolutions of the 
{\it set} $B_t$ of particles, rather than the {\it number} of  
particles. 
We shortly
describe the evolution, leaving a more formal definition for 
subsequent sections.
Let ${\bf q}=(q_t)_{t \in \N} $ be a sequence of probability measures 
on $\N$ such that
\bdnl{q_t}
q_t (1) <1\; \; \mbox{and}
\; \; \; m_t\st{\rm def.}{=}\sum_{k \ge 1}kq_t (k)<\8, 
\; \; \; \mbox{for all $t \in \N$}.
\edn
Then, the branching process with offspring
distribution ${\bf q}=(q_t)_{t\in \N}$, whose law is denoted by  
$P^{\bf q}$, is described
as the following dynamics:
\bitemize 
\item 
At time $t=0$, there is one particle ($|B_0|=1$).
\item 
Each particle in $B_t$ will die at time $t+1$, 
leaving a random number of children with law $q_t$ 
and all these variables are independent.
\eitemize
The {\it Galton-Watson process} can be thought of as the 
simplest process defined as above, 
where the offspring distribution is time-independent, $q_t \equiv q$ 
for all $t \in \N$. We then write $P=P^{\bf q}$ in short.
As is well known, the answer to the questions  
(\ref{GWsurv})--(\ref{GWgrow}) for the Galton-Watson process 
is given by:
\Theorem{GW}
\bds
\item[a)] \cite{AtNe72}
$P({\rm survival})>0$ if and only if 
$m\st{\rm def}{=}\sum_{k \ge 1}kq(k)>1$.
\item[b)] \cite{AtNe72,KeSti66a}
Suppose that $m >1$ and 
$
\sum_{k \ge 1}q(k) k \ln k  <\8.
$
Then
$$
\{ {\rm survival}\}=\Big\{ \lim_{t \ra \8}|B_t|/m^t \in (0,\8) \Big\},
\; \; \; \mbox{$P$-a.s.}
$$
\eds
\end{theorem}

On the other hand, 
the  {\it Smith-Wilkinson process} is defined as above, 
with ${\bf q}=(q_t)_{t \in \N} $  
a sequence of i.i.d. random probability measures 
on $\N$. We denote the law of ${\bf q}$ by $Q$ and set 
$$
P(\; \cdot\; )=\int Q (d{\bf q})P^{\bf q}(\; \cdot\; ).
$$
The answer to the questions  
(\ref{GWsurv})--(\ref{GWgrow}) for the Smith-Wilkinson process 
was given by W. L. Smith and W. Wilkinson \cite{SW69}, 
and by K. B. Athreya and S. Karlin \cite{AtKa71a,AtKa71b}.
\Theorem{SW}
Suppose for simplicity that 
$Q [|\ln m_t|]<\8$ (cf. (\ref{q_t})).
Then,
\bds
\item[a)]
\cite[Theorem 3.1]{SW69}, \cite[Theorems 1--3]{AtKa71a}. 
Either 
$$
\mbox{$P^{\bf q}({\rm survival})>0$, $Q$-a.s.}
\; \; \; \mbox{or}\; \; \; 
\mbox{$P^{\bf q}({\rm survival})=0$, $Q$-a.s.}
$$
The former case holds if and only if 
\bdnl{cond_SW}
Q [\ln m_t]>0\; \; \mbox{and}\; \; 
Q \ln {1 \over 1-q_t(0)}<\8.
\edn
\item[b)] \cite[Theorem 1]{AtKa71b}. 
In addition to (\ref{cond_SW}), suppose that
\bdnl{klnkSW}
Q \lef[m_t^{-1}\sum_{k \ge 1}q_t(k) k \ln k \ri] <\8.
\edn
Then, 
$$
\{ {\rm survival}\}
=\Big\{ \lim_{t \ra \8}{|B_t| \over m_0 \cdots m_{t-1}} \in (0,\8) \Big\},
\; \; \; \mbox{$P$-a.s.}
$$
\eds
\end{theorem}
\subsection{Branching random walks in random environment (BRWRE)}
 \label{sec:brw}
We now introduce the model of interest.  
To each $(t,x)
\in \N \times \zd$, we associate a distribution
$q_{t,x}=(q_{t,x}(k))_{k \in \N} \in \cP (\N)$ on the integers, which serves\
as an environment. Given ${\bf q}=(q_{t,x}; (t,x)
\in \N \times \zd)$, we define the branching random walk (BRW),
and we denote by  $P^{\bf q}$ its law, 
as the following dynamics:

\bitemize 
\item 
At time $t=0$, there is one particle 
at the origin $x=0$. 
\item 
Each particle, located at
site $x \in \zd$ at time $t$, jumps at time $t+1$ to one of the $2d$ neighbors 
of $x$ chosen uniformly; upon  arrival, it dies, leaving $k$ new particles 
there with probability 
$q_{t,x}(k)$. The number of newborn particles is independent of the jump,
and all these variables, indexed by the full population at time $t$, 
are independent. 
\eitemize

In the model of BRW with space-time random environment, 
we assume that ${\bf q}=(q_{t,x}; (t,x)
\in \N \times \zd)$ is an i.i.d. sequence in $\cP(\N)$ with some distribution
$Q$. We set
$P=\int Q (d{\bf q})P^{\bf q}$ as before. 

As is already mentioned, 
we denote by $B_t$ the set of particles present in the 
population at time $t$, and $B_{t,x}$ the set of 
those which are at site $x$ at time $t$.  
Though this model  is the same as in \cite{HuYo07,Yo08a},
we will formulate it in a detailed manner,  by representing 
each particle by its genealogy. This will provide some nice monotonicity
properties, like (\ref{branch<B}).

\subsection{The associated directed polymers in random environment} 
\label{sec:dpre}
From here on, we assume that 
\bdnl{hyp1}
Q[m_{0,0}+m_{0,0}^{-1}] < \8,
\; \; \mbox{where}\; \; 
m_{t,x}=\sum_{k \in \N}kq_{t,x}(k), \; \; (t,x) \in \zd.
\edn 
We write 
\bdnl{m}
m\st{\rm def.}{=}Q[m_{0,0}]<\8. 
\edn 
Let $(S_t)_{t \in \N}$ be the symmetric simple random walk on 
$\zd$ starting from $S_0 \equiv 0$ (cf. (\ref{1/2d})). 
We assume that $(S_t)_{t \in \N}$ is 
defined on a probability space $(\W_S, \cF_S, P_S)$. 
We then introduce the partition functions of 
{\it directed polymers in random environment} (DPRE),
\bdnl{Z}
Z_{t,x}=P_S[\z_t: S_t=x]
\; \; \; \mbox{and}\; \; \; 
 Z_t=P_S[\z_t],
\edn
where
\bdnl{z_t}
\z_t=\prod_{u=0}^{t-1}m_{u,S_u},
\edn
with $m_{t,x}$ from (\ref{hyp1}). The polymer measure is the law with weights 
$\z_t$ on the path space, and
a directed polymer in the environment ${\bf q}$ is a path sampled from this
law. 
It is easy to see (e.g., \cite[Lemma 1.3.1]{Yo08a}) that
\bdnl{DP=BRW}
Z_{t,x}=P^{\bf q}[|B_{t,x}|]\; \; \mbox{and}\; \;
 Z_t=P^{\bf q}[|B_t|].
\edn
For more information on DPRE, see \cite{CaHu03,CSY03,CSY04} for example. 
In standard literature on DPRE, 
the product in (\ref{z_t}) is defined as 
$\prod_{u=1}^{t}$ rather than $\prod_{u=0}^{t-1}$. 
This does not change large $t$ asymptotics, 
and the results we quote from there in the sequel are 
not affected by this slight difference.

\vvs
We set 
\bdnl{D}
\D =\{ \tht \in \R^d \; ; \; |\tht | \le 1 \}
\edn 
and,  for $\tht \in \D \cap \Q^d$,
\bdnl{m(tht)}
\N^* (\tht)=\{ t \in \N^* \; ; \; t\tht \in \zd, \; t-t|\tht| \in 2\N \},
\; \; \N (\tht)=\{0\} \cup \N^* (\tht).
\edn
We have $P_S (S_t=t\tht) \ge (2d)^{-t}>0$ 
for all $t \in \N (\tht)$. Note also that, 
\bdnl{n(tht)}
\N (\tht)=n(\tht )\N , \; \; \mbox{with $n(\tht ) =\min \N^* (\tht)$.}
\edn
Observe that (\ref{hyp1}), combined with
the elementary inequality 
$|\ln u| \le u \vee u^{-1}$ for $u >0$,
 implies that
\bdnl{lnZL^1}
\barray{ll}
Q |\ln Z_t| <\8 & \mbox{for all $t \ge 1$,} \\
Q |\ln Z_{t,x}| <\8 & \mbox{for all $(t,x) \in \N^* \times \zd$ 
with $P_S(S_t =x) >0$.}
\earray
\edn
\Proposition{press-cahu}
There exists a concave, 
upper semi-continuous function 
$
\psi: \D \ra \R
$, with $\Delta$ from (\ref{D}), such that 
\bdnl{psi(tht)} 
\psi (\tht )=\lim_{t \ra \8 \atop t \in \N^* (\tht)} 
\frac{1}{t} Q[\ln Z_{t,t\tht }]
\; \; \mbox{for all $\tht \in \D \cap \Q^d$.}
\edn 
$\psi$ is symmetric in the sense that 
$\psi (\tht)=\psi (|\tht_{\s(1)}|, \ldots,|\tht_{\s(d)}|)$ for any 
permutation $\s$ of $\{ 1,\ldots,d\}$, and
\bdnl{psi(e)psi(0)}
Q[\ln m_{0,0}]-\ln (2d) =\psi (\pm e_i) 
\le \psi ( \tht ) \le \psi ( 0 ).
\edn
 Moreover, the following limit 
also exists and equals $\psi (0)$:
\bdnl{psi(0)} 
\Psi \st{\rm def.}{=}\lim_{t \ra \8}\frac{1}{t} Q[\ln Z_t].
\edn
\end{proposition}
We refer the reader to \cite[page 287, Theorems 1.1 and 1.2]{CaHu03} 
for the proof of \Prop{press-cahu} (Although the 
random variable $\ln m_{0,0}$ is assumed to be Gaussian 
in \cite{CaHu03}, it is not essential for this purpose).
The number $\Psi$ is called the free energy of the polymer, and
$\psi ( \tht )$  the directional free energy.

\SSC{The results}
\subsection{Criteria for global and local survival, 
and growth rates} \label{sec:main}
We suppose (\ref{hyp1}).
Then, the growth rate 
for the total population 
of the BRWRE is identified with the free energy (\ref{psi(0)})
of the DPRE associated to it by (\ref{Z}). 
\Theorem{N_T}{\bf(global growth)}
Let $\e>0$.
\bds
\item[a)] 
We have, $P$-a.s.,
\bdnl{N_T<1}
|B_t|  \le e^{(\Psi+\e)t}\; \; \mbox{for all large $t$'s.}
\edn
More precisely, there exist $c_1,c_2  \in (0,\8)$ such that 
\bdnl{N_T<2}
P \lef( |B_t| \ge e^{(\Psi+\e)t} \ri) \le 3 \exp (-c_1\e^2t)\; \; \; 
\mbox{for all $\e \in (0,c_2]$, $t \ge 1$.}
\edn
In particular,
$$
P( \mbox{\rm survival} )=0 \; \; \mbox{if}\; \; \Psi <0,
\; \; \; \mbox{cf. (\ref{DefSurv}).}
$$
\item[b)]
Suppose $\Psi >0$ and 
\bdnl{1-q(0)}
Q\ln {1 \over 1-q_{0,0}(0)} <\8.
\edn
Then, 
\bdnl{surv}
P^{\bf q}( \mbox{\rm survival} ) >0,\; \; \mbox{Q-a.s.}
\edn
Suppose $\Psi >0$, (\ref{1-q(0)}) and 
\bdnl{KS}
Q \lef[ {P^{\bf q}[|B_{t,0}|\ln |B_{t,0}|]
 \over Z_{t,0}}\ri]<\8
\; \; \mbox{for all $t \in 2\N^*$.}
\edn
Then, 
\bdnl{N_T>}
\{ \mbox{\rm survival} \}=
\Big\{ |B_t|  \ge e^{(\Psi-\e)t}\; \; \mbox{for all large $t$'s} \Big\},
\; \; \; \mbox{$P$-a.s.}
\edn
\eds
\end{theorem}
In particular, (\ref{N_T>}), together with (\ref{N_T<1}) says that 
$\Psi >0$ implies 
\bdnl{survie-log}
\{ \mbox{\rm survival} \}=
\Big\{ \lim_{t \ra \8}{1 \over t}\ln |B_t|=\Psi\;  \Big\},
\; \; \; \mbox{$P$-a.s.}
\edn
under the conditions (\ref{1-q(0)}) and (\ref{KS}). 
The assumption ``$\Psi >0$ and (\ref{1-q(0)})" for \Thm{N_T}(b) is the 
generalization of (\ref{cond_SW}) (also of ``$m>1$" in \Thm{GW}), 
while the condition 
(\ref{KS}) is the counterpart of (\ref{klnkSW}). 
We will explain in the remarks 1)--2) below 
that the assumptions in \Thm{N_T}(b) are not too restrictive, 
while, in remark 3), 
 we compare the survival probability $P({\rm survival})$ 
with those of the 
Galton-Watson process and the Smith-Wilkinson process, 
which are properly associated with our BRWRE.

\vvs 
\noindent {\bf Remarks:} 
{\bf 1)} 
Assume that $m_{0,0}$ is not a constant a.s. Then, 
by \Prop{press} (b), 
$Q[\ln m_{0,0}] \ge 0$ (for which $m_{0,0} \ge 1$ a.s. is more than 
enough) implies that $\Psi >0$. 

\vs
\noindent {\bf 2)} 
We will see at the end of section \ref{sec:main} that 
the following mild integrability condition implies 
both (\ref{1-q(0)}) and (\ref{KS}).
\bdnl{K^2}
Q\lef[\ln m^{(2)}_{0,0}\ri]<\8, 
\; \; \mbox{where}\; \; m^{(2)}_{t,x}=\sum_{k \in \N}k^2q_{t,x}(k).
\edn  

\vs
\noindent {\bf 3)} 
Let us suppose (\ref{1-q(0)}) and 
let $\s^{\rm GW}$ and $\s^{\rm SW}$ be 
survival probabilities 
for Galton-Watson model with offspring distribution 
$Q[q_{0,0}(\; \cdot \; )]$ and the 
Smith-Wilkinson process with the 
offspring distributions $q_{t,0}$, $t \in \N^*$.
By computations of 
generating function in section \ref{gf} below, we will see that 
\bdnl{SW<surv<GW}
\s^{\rm SW} \le P({\rm survival}) \le \s^{\rm GW}.
\edn
On the other hand, it is known that 
$\s^{\rm SW}>0$ if and only if $Q [\ln m_{0,0}]>0$ \cite{AtKa71a,SW69}. 
Thus, by \Thm{N_T}, we have 
$$
\s^{\rm SW} =0 < P({\rm survival})
\; \; \mbox{if $Q [\ln m_{0,0}] \le 0 <\Psi$.}
$$ 
Similarly, we have 
$$
P({\rm survival})=0 <  \s^{\rm GW}\; \; \mbox{if $\Psi <0<\ln m$}.
$$
We now turn to the criterion for the local survival and the identification 
of the local growth rate in terms 
of the directional growth rate (\ref{psi(tht)}): 
\Theorem{N_Tloc}{\bf(local growth)}
Let $\tht \in \D \cap \Q^d$ and $\e>0$.
Then, 
\bds
\item[a)]
$P$-almost surely,
\bdnl{N_T<loc1}
|B_{t,t\tht}|  \le e^{(\psi(\tht)+\e)t}\; \; 
\mbox{for all large $t \in \N $.}
\edn
More precisely, there exist $c_1,c_2  \in (0,\8)$ such that 
\bdnl{N_Tloc<2}
P \lef( |B_{t,t\tht}| 
\ge e^{(\psi(\tht)+\e)t} \ri) \le 3 \exp (-c_1\e^2t)
\edn
for all $\e \in (0,c_2]$ and $t \in \N $.
In particular,
$$
P( B_{t,t\tht} \neq \epty \; \;  \mbox{i.o.} )=0
 \; \; \mbox{if}\; \; \psi (\tht ) <0,
$$
\item[b)]
Supppose $\psi(\tht ) >0$ and (\ref{1-q(0)}). Then, 
\bdnl{loc_surv9}
P^{\bf q}\big( B_{t,t\tht} \neq \epty \; \; \mbox{for all 
$t \in \N (\tht)$} \big)>0,
\; \; \mbox{$Q$-a.s.}
\edn
Supppose $\psi(\tht ) >0$, (\ref{1-q(0)}) and
\bdnl{KSloc9}
Q \lef[ {P^{\bf q}[|B_{t,t\tht}|\ln |B_{t,t\tht}|]
\over Z_{t,t\tht}}\ri]<\8
\; \; \mbox{for all $t \in \N^*(\tht )$.}
\edn 
Then, 
\bdnl{N_T>loc9}
\{ B_{t,t\tht} \neq \epty \; \; \mbox{for all 
$t \in \N (\tht)$}\} \sub 
\{ |B_{t,t\tht}|  \ge e^{(\psi(\tht )-\e)t}\; \; 
\mbox{for all large $t \in \N (\tht)$} \}, 
\; \; \mbox{$P$-a.s.}
\edn
\eds
\end{theorem}
\noindent {\bf Remarks:} 
{\bf 1)} 
Since $\psi :\D \ra \R$ is concave, 
it is continuous in the interior of $\D$ \cite[page 82, Theorem 10.1]{Roc70}. 
Hence, if $\Psi =\psi (0)>0$, then $\psi (\tht )>0$ in a 
neighborhood of 0. 

\vs
{\bf 2)} By (\ref{N_T<loc1}) and (\ref{N_T>loc9}) the local growth rate is
given by $\psi(\theta)$ when it is positive.
We will see in Remark 1) below  \Prop{press} that $\psi(\theta)$
is equal to the global growth rate minus the large deviations 
rate function at $\theta$ for the polymer measure. This relation  
agrees with, and extends, the one for standard branching random walks as well
as Example 1.10 in \cite{CP07} for constant branching.

\vs
{\bf 3)}
The mild integrability condition
(\ref{K^2}) implies (\ref{1-q(0)}) and (\ref{KSloc9}) for 
all $\tht \in \D \cap \Q^d$. We now prove these claims.

Note that $1-q_{0,0}(0)=\sum_{k \ge 1}q_{0,0}(k)$ 
and hence that
$$
m_{0,0}^2 = \lef( \sum_{k \ge 1}kq_{0,0}(k) \ri)^2 
\st{\mbox{\scriptsize Schwarz}}{\le} (1-q_{0,0}(0))m^{(2)}_{0,0}.
$$
Therefore,
$$
\ln {1 \over 1-q_{0,0}(0)} 
\le \ln {m^{(2)}_{0,0} \over m_{0,0}^2}
= \ln m^{(2)}_{0,0}-2\ln m_{0,0},
$$
which proves that (\ref{K^2}) implies (\ref{1-q(0)}).

To show (\ref{KSloc9}), we note that
$$
{P^{\bf q}[|B_{t,t\tht}|\ln |B_{t,t\tht}|]
\over Z_{t,t\tht}}
\st{\mbox{\scriptsize Jensen}}{\le}
\ln {P^{\bf q}[|B_{t,t\tht}|^2] \over Z_{t,t\tht}} 
=\ln P^{\bf q}[|B_{t,t\tht}|^2]-\ln Z_{t,t\tht}.
$$
Thus, in view of (\ref{lnZL^1}), it is enough to prove that
\bds
\item[1)] \hspace{1cm}
${\dps Q\ln P^{\bf q}[|B_t|^2] <\8}$ for all $t \ge 1$.
\eds
To this end, we write
$$
|B_t|^2 =  \lef( \sum_{(x,\n) \in B_{t-1}}K_{t-1,x,\n} \ri)^2
\le |B_{t-1}|\sum_{(x,\n) \in B_{t-1}}K_{t-1,x,\n}^2,
$$
and 
$$
P^{\bf q}[|B_t|^2 |\cF_{t-1}] 
\le |B_{t-1}|\sum_{(x,\n) \in B_{t-1}}m_{t-1,x}^{(2)}
\le |B_{t-1}|^2\max_{|x| \le t-1}m_{t-1,x}^{(2)}.
$$
By iteration, we get
$$
P^{\bf q}[|B_t|^2 ] 
\le \prod_{u=0}^{t-1}
\max_{|x| \le u}m_{u,x}^{(2)},
$$
and hence
$$
\ln P^{\bf q}[|B_t|^2]
\le \sum_{u=0}^{t-1}
\max_{|x| \le u}\ln m_{u,x}^{(2)}
\le \sum_{u=0}^{t-1}
\sum_{|x| \le u}\ln m_{u,x}^{(2)},
$$
which proves that (\ref{K^2}) implies 1).
\hfill $\Box$

\subsection{More on the survival probability} 
There are still additional observations which we can make on the 
survival probability. 
The first of these is the following zero-one law:
\Proposition{surv01}
Either
$$
P^{\bf q}({\rm survival})>0\; \; \mbox{$Q$-a.s.}
\; \; \; \mbox{or}\; \; \; 
P^{\bf q}({\rm survival})=0\; \; \mbox{$Q$-a.s.}
$$
\end{proposition}
The proof of \Prop{surv01} is presented in section \ref{gf}.
We already know the zero-one law when $\Psi \neq 0$ from Theorem \ref{Thm.N_T},
but the result is new in the critical case $\Psi=0$. 
Concerning the survival/extinction at criticality we leave the 
following conjecture for future work:
\begin{conjecture}
When $\Psi=0$, $P({\rm survival})=0$.
\end{conjecture}
The conjecture seems quite plausible, since the 
genealogy in BRWRE is phenomenologically 
quite similar to the active paths in the contact process 
and open oriented paths in the oriented percolation. For 
the contact process and the oriented percolation, 
the equivalent of the above conjecture is known to be true 
\cite{BG90,GH02}. 

Another observation is concerned with 
the non-degeneracy of the limit of a martingale, 
 which is obtained by normalizing the total population $|B_t|$. 
To state it, we set 
\bdnl{ovn} 
W_t=|B_t|/m^t, 
\edn 
where $m$ is defined by (\ref{m}).
The non-negative sequence $W_t$ is a 
$(P, (\cF_t))$-martingale, and
by Doob's martingale convergence theorem \cite{Dur95},
the following limit exists:
\bdnl{ovn_8}
W_\8=\lim_t W_t, \; \; \mbox{$P$-a.s.} 
\edn 
It is known \cite[Corollaries 1.2.2 and 3.3.2]{Yo08a} that
$$
P(W_\8 >0)\lef\{\barray{ll}
>0 & \mbox{if (\ref{WD}), $m >1$ and $Q[m^{(2)}_{0,0}]<\8$}, \\
=0 & \mbox{if (\ref{SD})}.
\earray \rig.
$$
In particular, there is a whole range of parameters where
$P(W_\8 >0)=0<P(\mbox{\rm survival})$, namely, when both
(\ref{SD}) and $\Psi>0$ hold. However,
arguments used to prove \Thm{N_T} (\Lem{701} below) also lead  to:

\Proposition{reg=surv}
Suppose that $m<\8$ and $P(W_\8 >0)>0$. Then,
$$
\{ \mbox{\rm survival} \}\st{\rm a.s.}{=}\{ W_\8>0\}
\st{\rm a.s.}{=}\Big\{ \lim_{t \to \8} \frac{|B_t|}{Z_t} \in (0, \8)
\Big\}.
$$
\end{proposition}

The above proposition is an analogue of a classical 
result for the Galton-Watson process 
\cite[page 9, Theorem 9 (iii)]{AtNe72} 
(Note that the hypothesis $\s^2 <\8$ is not used there, but only $m < \8$).
It extends to BRWRE, \Thm{GW} (b) and  
\Thm{SW} (b). Under an additional assumption, it improves 
(\ref{survie-log}).
We also remark that the equality analogous to \Prop{reg=surv} 
holds true for continuous and discrete time 
linear systems \cite{Gri83,Yo08b} taking values in $\N^{\zd}$.
Finally, note that non-degenerate limits for local ratios
$|B_{t,x}|/Z_{t,x}$ have been found for homogeneous
BRW \cite{Big79}, but in the case of random environment 
this seems to be beyond reach.

\SSC{Preliminaries}
\subsection{Construction of BRWRE}
To give a precise definition of the above dynamics, we will adopt here 
a reformulation given in \cite{Na09}, and introduce a tree 
which set of vertices consists of all possible ancestral histories of 
the branching process: 
\bdnl{cT}
\cT=\bigcup_{t \in \N}\cT_t
\; \; \; 
\mbox{where $\cT_t =\{ \n =(\n_i)_{i=0}^t \in (\N^*)^{1+t}\; ; \; \n_0=1\}$.}
\edn
For $\n \in \cT$, $|\n| \in \N$ stands for its generation 
defined by $\n \in \cT_{|\n|}$. 
For $\n \in \cT_{t+s}$, $\n |_t$ denotes the 
ancestral history of $\n$ up to generation $t$:
\bdnl{n|_t}
\n |_t =(1,\n_1,\ldots,\n_t) \in \cT_t.
\edn

Let $p(\cdot, \; \cdot)$ be the transition probability for
the symmetric simple random walk on $\zd$: 
\bdnl{1/2d}
p(x,y)=\left\{ \barray{ll}
\frac{1}{2d} & \mbox{if $|x-y|=1$,} \\
0 & \mbox{if $|x-y|\neq 1$.}\earray \rig. \edn 
Fix an environment ${\bf q}=(q_{t,x})_{(t,x) \in \N \times \zd}$,
$q_{t,x} \in \cP (\N)$.
 Then, the branching random walk (BRW) with offspring
distribution ${\bf q}$ is described in an unformal manner
as follows:

\bitemize 
\item 
At time $t=0$, there is one particle 
at the origin $x=0$ (the founding ancestor labeled by "1"). 
\item 
Suppose that there is a particle at
site $x \in \zd$ at time $t$, which is 
labeled by $\n=(1,\n_1,..,\n_t)$, $\nu_i \in \N^*$. 
At time $t+1$, it jumps to a site $y$
with probability $p(x,y)$ independently of the other particles.
 At arrival, it dies, leaving $k$ new particles there with probability 
$q_{t,x}(k)$. The newborn particles are 
labeled by $(1,\n_1,\ldots,\n_t,\n_{t+1})$ 
(with $\n_{t+1}=1,\ldots,k$).
\eitemize

To put it on steady grounds, 
we introduce a few ingredients.

\vvs \noindent $\bullet$ {\it Spatial motion:}
We define the measurable
space $(\W_X,\cF_X)$ as the set $(\zd)^{\N \times \zd \times \cT}$ 
with the product $\s$-field, and 
$\W_X \ni X \mapsto X_{t,x,\n}$ 
for each $(t,x,\n) \in  \N \times \zd \times \cT$ as
the projection. We define $P_X \in \cP (\W_X,\cF_X)$ as the
product measure such that 
\bdnl{P_X} 
P_X (X_{t,x,\n}=y)=p(x,y)\;
\; \; \mbox{for all $(t,x,\n) \in \N \times \zd \times \cT$ and
$y \in \zd$.} \edn 
Here, we interpret $X_{t,x,\n}$ as the position
at time $t+1$ of the children born from the particle 
labeled by $\n$, when it occupies the time-space location $(t,x)$.

\vvs 
\noindent $\bullet$ {\it Offspring distribution:} 
We set
$\W_q =\cP (\N)^{\N \times \zd}$.
Thus, each ${\bf q} \in \W_{\bf q}$ is a function $(t,x) \mapsto
q_{t,x}=(q_{t,x}(k))_{k \in \N}$ from $\N \times \zd$ to $\cP
(\N)$. We interpret $q_{t,x}$ as the offspring distribution for
each particle which occupies the time-space location $(t,x)$. The
set $\cP(\N)$ is equipped with the natural Borel $\s$-field
induced from that of $[0,1]^\N$. We denote by $\cF_{\bf q}$ the product
$\s$-field on $\W_{\bf q}$.

We define the measurable space $(\W_K,\cF_K)$ 
as the set $\N^{\N \times \zd \times \cT}$ 
with the product $\s$-field, and $\W_K
\ni K \mapsto K_{t,x,\n}$ for each $(t,x,\n) \in \N \times \zd \times \cT$ 
as the projection. For each fixed ${\bf q} \in \W_{\bf q}$, we
define $P^{\bf q}_K \in \cP (\W_K,\cF_K)$ as the product measure such
that 
\bdnl{P^q_K} 
P^{\bf q}_K (K_{t,x,\n}=k)=q_{t,x}(k)\; \; \;
\mbox{for all $(x,t,\n) \in \N \times \zd \times \cT$ and $k \in
\N$.} \edn 
We interpret $K_{t,x,\n}$ as the number of the children
born from the particle $\n$ at time-space location $(t,x)$.

\vvs 
\noindent 
$\bullet$ {\it Branching random walk in random environment} (BRWRE): 
We fix a product measure $Q \in \cP (\W_{\bf q},\cF_{\bf q})$, which describes
the i.i.d. offspring distribution assigned to each time-space
location. 
Finally, we define $(\W, \cF)$ by
$$
\W=\W_X \times \W_K \times \W_{\bf q}, \; \; \; \cF=\cF_X \otimes \cF_K
\otimes \cF_{\bf q},
$$
and $P^{\bf q}, P \in \cP (\W, \cF)$ by
$$
P^{\bf q}=P_X \otimes P^{\bf q}_K \otimes \del_{\bf q}, 
\; \; \; P=\int Q (d{\bf q})P^{\bf q}.
$$
We define a Markov chain $(B_t)_{t \in \N}$ with 
values in finite subsets of $\zd \times \cT$, inductively by 
$B_0=(0,1)$, and, for $t \ge 1$, 
 \bdnl{B_t}
B_t=\bigcup_{(x,\n) \in B_{t-1}}
\Big\{ (y,\m)  \in \zd \times \cT_t\; ; \; 
 X_{t-1,x, \n}=y, \; \m |_{t-1}=\n, \;  \m_t \le K_{t-1,x, \n} \Big\}
.
 \edn
We call the above process 
{\it the branching random walk in random environment} (BRWRE). 
It is meant by $(x,\n) \in B_t$ that the time-space location 
$(t,x)$ is occupied by a particle with its ancestoral history $\n$.  
We consider the filtration: 
\bdnl{cF_t} \cF_0=\{ \epty, \W\},\; \;
\cF_t =\s (X_{s, \cdot,\cdot}, K_{s, \cdot,\cdot}, q_{s, \cdot}\;
; \; s \le t-1 ),\; \; \; t \ge 1, \edn 
which the process $(B_t)_{t \in \N}$ is adapted to. 
We define, for fixed $(t,x)\in \N \times \zd$, 
the set of particles in $B_t$, which 
occupies the site time-space $(t,x)$ is denoted by:
\bdnl{N_t}
 B_{t,x}=\{ (y,\n) \in B_t\; ; \; y=x\}.
\edn
We remark that the total population $|B_t|$ is exactly the classical
Galton-Watson process if $q_{t,x} \equiv q$, where $q \in \cP
(\N)$ is non-random. On the other hand, if $\zd$ is replaced a
singleton, then $|B_t|$ is the population of the Smith-Wilkinson
model \cite{SW69}.

\vvs
\noindent {\bf Remark:}
The definition (\ref{N_t}) is consistent with that 
in \cite{HuYo07,Yo08a}. In fact, it is easy 
to see from (\ref{B_t}) that 
\bdnl{indN_t}
|B_{0,y}|=\del_{0,y},
\; \; |B_{t,y}|=\sum_{(x,\n ) \in B_{t-1}}
 \del_y (X_{t-1,x,\n})K_{t-1,x,\n}, \; \; t \ge 1.
\edn
If we write 
$$
\{ \n (i) \}_{i=1}^{|B_{t-1,x}|}
=\{ \n \in \cT_{t-1} : (x,\n ) \in B_{t-1}\}
$$
for each $x \in \zd$, then, (\ref{indN_t}) becomes:
\bdnl{indN_t2}
|B_{0,y}|=\del_{0,y},
\; \; |B_{t,y}|=\sum_{x \in \zd}\sum_{i=1}^{|B_{t-1,x}|}
 \del_y (X_{t-1,x,\n (i)})K_{t-1,x,\n (i)}, \; \; t \ge 1,
\edn
which is the recursion used in \cite{HuYo07,Yo08a}.
\subsection{Complements on DPRE}
We denote by $I_S$ the large deviation rate 
function for the random variables $(S_t/t)_{t \in \N^*}$:
\bdnl{I_S}
I_S (\tht)=\sup_{\a \in \rd}
\{ \a \cdot \tht -\ln P_S[\exp (\a \cdot S_1)]  \},
\; \; \; \tht \in \D.
\edn
It is well-known that
\bdnl{<I_S<}
0=I_S (0) \le I_S (\tht) \le I_S (\pm e_i)=\ln (2d),
\edn
where $e_i=(\del_{i,j})^d_{j=1}$,
and,
for all $\tht \in \D \cap \Q^d$, that 
\bdnl{-I_S}
-I_S (\tht)=\lim_{t \ra \8 \atop t \in \N^*(\tht)}
\frac{1}{t}\ln P_S(S_t =t\tht)
=\sup_{t \in \N^*(\tht)}\frac{1}{t}\ln P_S(S_t =t\tht)
.
\edn
\Proposition{press}
\bds
\item[a)] 
We have
\bdnl{<psi<}
Q [\ln m_{0,0}] -I_S (\tht) \le \psi (\tht) \le \ln m -I_S (\tht)
\; \; \mbox{for all $\tht \in \D \cap \Q^d$.}
\edn
(Note that $Q[|\ln m_{0,0}|] <\8$ by (\ref{hyp1})).
 Moreover, the first inequality 
is strict if $\tht =0$ and 
$m_{0,0} \not\equiv m$ a.s. Finally, the first inequality in (\ref{<psi<})
 is an equality if 
$\theta= \pm e_i$.
\item[b)]
There are constants $c_1,c_2\in (0,\8)$ such that 
\bdnl{exp_con}
Q \lef(  {1 \over t}\lef| \ln Z_t  
-Q \lef[ \ln Z_t\rig]\rig| >\e \ri)
  \le  2\exp \lef( -c_1\e^2t \rig),
\edn
for all $\e \in (0,c_2]$ and $t \in \N^*$, and 
\bdnl{exp_con_loc}
Q \lef(  {1 \over t}\lef| \ln Z_{t,t\tht }  
-Q \lef[ \ln Z_{t,t\tht }\rig]\rig| >\e \ri)
 \le  2\exp \lef( -c_1\e^2t \rig),  
\edn 
for all $\e \in (0,c_2]$, $\tht \in \D \cap \Q^d$ and $t \in \N^* (\tht)$.
As a consequence, 
\bdnl{psi_a.s.}
\lim_{t \ra \8} \frac{1}{t} \ln Z_t=\Psi 
\; \; \mbox{and}\; \; 
\lim_{t \ra \8 \atop t \in \N^* (\tht)} 
\frac{1}{t} \ln Z_{t,t\tht }=\psi (\tht )
,\; \; 
\mbox{$Q$-a.s.}
\edn
\eds
\end{proposition}
We prove \Prop{press}(a) in section \ref{ppressb}.
A way to prove \Prop{press}(b) is to apply a series of general martingale 
inequalities as is done in \cite{LW08}, theorem 7.2. For the reader's 
convenience, we enclose a short proof in section \ref{ppressc}. 

\vvs
\noindent {\bf Remarks:}
{\bf 1)} By (\ref{psi_a.s.}) (also as is discussed in \cite{CaHu03}), 
the function $\psi (0)-\psi (\cdot) : \D \ra [0,\8)$ gives,
for almost all realizations of the environment, 
the large deviation rate function for the random probability 
measures (usually called polymer measures):
$$
{1 \over Z_t}P_S\lef[\z_t : S_t/t \in \cdot \ri]
\; \; \; \mbox{as $t \ra \8$}.
$$
\noindent {\bf 2)}
The quantity $\Psi$ is 
called the free energy in the context of DPRE. 
Its value relative to $\ln m$ is important there, and hence is  
well studied. $\Psi \le \ln m$ by \Prop{press}(b).  
$\Psi=\ln m$ if
\bdnl{WD}
\mbox{$d \ge 3$ and 
${\dps {Q[ m_{0,0}^2] \over m^2} <{1 \over \pi_d}}$},
\edn
where $\pi_d=P_S (S_t =0 \; \mbox{for some $t \ge 1$})$. 
More precisely, $(Z_t/m^t)_{t \ge 1}$ converges to a 
positive limit $Q$-a.s. under the condition (\ref{WD})
(\cite[page 128, Remark 3.2.3]{CSY04}, 
\cite[page 282, Lemma 1]{SoZh96}).
On the other hand, $\Psi < \ln m$ if 
\bdnl{SD}
\barray{l}
\mbox{$d=1,2$ and $m_{0,0}$ is not a constant,}\\
\mbox{or}\; \; 
{\dps Q\lef[ {m_{0,0} \over m}\ln {m_{0,0} \over m} \ri]>\ln (2d)}
\earray
\edn
(\cite[page 709, Theorem 2.3(a)]{CSY03}, \cite{CV04,La09}).

\SSC{Proofs}

\subsection{Proof of \Thm{N_T}(a) and \Thm{N_Tloc}(a)}
We first show (\ref{N_T<2}). 
By superadditivity (e.g.,\cite[page 720, Proof of Proposition 2.5]{CSY03}) 
and \Prop{press}, we have for large $T$ that
$$
Q[\ln Z_T] \ge T \Psi
\; \; \mbox{and}\; \; 
Q \lef( \ln Z_T \ge Q[\ln Z_T] +{\e \over 2}T \ri) \le 2 \exp (-c_1\e^2T).
$$
Thus, 
\bds
\item[1)]
$Q \lef( Z_T \ge e^{(\Psi+{\e \over 2})T} \ri) 
\le 2 \exp (-c_1\e^2T).$
\eds 

On the other hand, we have  
for any $\chi =\chi ({\bf q}) \in \{0,1\}$ that
\bdnn
P^{\bf q} \lef( |B_T| \ge e^{(\Psi+\e)T} \ri) 
& \le &
P^{\bf q} \lef( |B_T| \ge e^{(\Psi+\e)T} \ri) \chi +1-\chi \\
& \st{\mbox{\scriptsize Chebyshev}}{\le} &
e^{-(\Psi+\e)T}Z_T\chi +1-\chi.
\ednn 
We choose $\chi ={\bf 1}\{ Z_T < e^{(\Psi+{\e \over 2})T}\}$ and take 
$Q$-expectation to get 
$$
P \lef( |B_T| \ge e^{(\Psi+\e)T} \ri) \le 
e^{-{\e \over 2}T}+Q \lef( Z_T \ge e^{(\Psi+{\e \over 2})T} \ri) 
\st{\scriptstyle 1)}{\le} e^{-{\e \over 2}T}+2 \exp (-c_1\e^2T).
$$
We now obtain (\ref{N_T<1}) by Borel-Cantelli lemma. 
This proves \Thm{N_T}(a). The proof of \Thm{N_Tloc}(a) is similar.
\hfill $\Box$
\subsection{Local survival and growth rate along subsequences of times}
The purpose of this subsection is to prove the following 
two lemmas, which proves (\ref{surv}) 
and will also be our first step for the proof of \Thm{N_Tloc}(b).
\Lemma{P>0}
Suppose $\tht \in \D \cap \Q^d$, (\ref{1-q(0)}) and that 
there exists $T \in \N^* (\tht )$ such that 
\bdnl{assP>0}
 Q[\ln Z_{T,T\tht}] >0.
\edn
Then, 
\bdnl{survP>0}
P^{\bf q}(B_{sT,sT\tht}\neq \epty \; \mbox{for all $s \ge 1$})>0,
\; \; \mbox{Q-a.s.}
\edn
Suppose in addition that
\bdnl{KS_P>0}
Q \lef[ {P^{\bf q}[|B_{T,T\tht}|\ln |B_{T,T\tht}|]
\over Z_{T,T\tht}}\ri]<\8.
\edn
Then, for any $1 \le r <\exp \lef( {1 \over T}Q[\ln Z_{T,T\tht}] \ri)$, 
\bdnl{claimP>0}
P^{\bf q} \lef( \inflim_{s \ra \8}r^{-sT}|B_{sT,sT\tht}| >0 \ri)>0,
\quad \mbox{$Q$-a.s.}
\edn
\end{lemma}
\Lemma{checkP>0}
Suppose $\tht \in \D \cap \Q^d$, $\psi (\tht) >0$, and (\ref{1-q(0)}).
Then (\ref{survP>0}) holds for 
any large enough $T \in \N^*(\tht)$. 
If we suppose (\ref{KSloc9}) in addition,  
then, (\ref{claimP>0}) holds for 
any large enough $T \in \N^*(\tht)$ and for any 
$r<\exp (\psi (\tht))$.
\end{lemma}
Since $\Psi=\psi (0)$, we get (\ref{surv}) from \Lem{checkP>0}.

\vvs
Let us first prove \Lem{checkP>0}, assuming \Lem{P>0}. \\
Let $1 \le r  <\exp (\psi (\tht))$. 
Then, by  (\ref{psi(tht)}),
$$
{1 \over T}Q [\ln  Z_{T, T\tht }]
> \ln r\; \; 
 \mbox{for all large $T \in \N^* (\tht)$,}
$$ 
and the claim follows from \Lem{P>0}. 
\hfill $\Box$

\vvs  
We now turn to the proof of \Lem{P>0}. The strategy is to 
imbed a Smith-Wilkinson process into 
$(B_{sT, sT\tht})_{s \in \N}$. 
This will 
reduce the study of survival of $(B_{sT, sT\tht})_{s \in \N}$ to that 
of the Smith-Wilkinson process, for which the criterion is known by 
\Thm{SW}. In fact, thanks to \Thm{SW}, the imbeded Smith-Wilkinson process 
can be made supercritical for large $T$.

To write down the above mentioned imbedding precisely, 
we introduce some notation. 
For $(s,z,\lm ) \in \N \times \zd \times \cT_s$, we define the 
{\it $(s,z,\lm)$-branch} 
$$
(B^{s,z,\lm}_t)_{t \in \N}
$$
of $(B_t)_{t \in \N}$  inductively by $B^{s,z,\lm}_0=(z,\lm)$ and 
for $t \ge 1$,
\bdnl{branch_def}
B^{s,z,\lm}_t=\bigcup_{(x,\n) \in B^{s,z,\lm}_{t-1}}
\{ (y,\m) \in \zd \times \cT_{s+t}  \; ; \; 
 X_{s+t-1,x, \n}=y, \; \m |_{s+t-1}=\n, \;  \m_{s+t} \le K_{s+t-1,x, \n} \}.
\edn
This amounts to restarting a BRWRE from a single 
particle at time-space $(s,z)$, whose 
ancestoral history up to time $s$ is given by $\lm$. Clearly,
\bdnl{branch<B}
\{ (z,\lm) \in B_s\} \sub \{ B^{s,z,\lm}_t \sub B_{s+t}\; 
\mbox{for all $t \ge 0$}\}.
\edn
Note also that $B^{s,z,\lm}_\cdot$ is a function of 
$$
\{ (X_{s+\cdot,\cdot, \n}, K_{s+\cdot,\cdot, \n}) \;; \; \n|_s=\lm \}.
$$
and therefore,  
for each fixed $(s,z)$, 
\bdnl{branch_ind}
\mbox{$\{ B^{s,z,\lm}_\cdot \}_{\lm \in \cT_s}$ are i.i.d. under 
$P^{\bf q}$.}
\edn
For $(t,x) \in  \N \times \zd$, the set of particles 
in $B^{s,z,\lm}_t$, which occupy the site $x$ is denoted 
by:
\bdnl{branch_N_t}
B^{s,z,\lm}_{t,x}
=\{ (y,\n) \in B^{s,z,\lm}_t\; ; \; y=x\}.
\edn
We fix a $T$ such that (\ref{assP>0}) holds. We will then define 
a Markov chain $(B^*_s)_{s \in \N}$ with values in finite 
subsets of $\cT$, which serves as a ``lower bound" 
of $(B_{sT, sT\tht})_{s \in \N}$. 
We now define $B^*_0=\{1\}$, and for $s \geq 1$,
\bdnl{eq:defN*}
 B^*_s =\bigcup_{\lm \in B^*_{s-1}}B^*_{s,\lm}, 
\; \; \mbox{with}\; \; 
B^*_{s,\lm}=\lef\{ \n \in \cT_{sT}
\; ; \; (sT\tht,\n) \in B^{(s-1)T,(s-1)T\tht,\lm}_T \ri\} .
\edn 
In words,  $B^*_s$ is the subset of $B_{sT, sT\tht}$ composed of genealogies
which have been at site $rT\tht$ at times $rT$ for $r=1,2,\ldots,s$.
\medskip

\Lemma{iid}
\bds
\item[a)] 
For each $s \ge 1$, $\{ B^*_{s,\lm} \}_{\lm \in \cT_{(s-1)T}}$ 
are i.i.d. under $P^{\bf q}$.
\item[b)] 
Define $q^*_s \in \cP (\N)$, $s \in \N$ by
\bdnl{eq:defq*}
q^*_s(k)=P^{\bf q}\lef( |B^*_{s,\lm}|=k\ri),
\; \; \; k \in \N.
\edn 
Then, the  sequence $(q^*_s)_{s \in \N}$ is i.i.d. under $P$. 
\eds
\end{lemma}
Proof: a) This follows from (\ref{branch_ind}).\\
b) Since $q^*_s$  is $\cF_{sT}$-measurable, it is enough to show that 
$$
P(q^*_s \in \cdot |\cF_{(s-1)T})=P(q^*_1 \in \cdot ).
$$
Let $\tht_{t,x} :\w \mapsto \tht_{t,x}\w$ be 
time-space shift for $\w=(X,K,q)$. Then, 
$q^*_s =q^*_1 \circ \tht_{(s-1)T,(s-1)T\tht}$. 
By the shift-invariance, we have 
$$
P(q^*_s \in \cdot |\cF_{(s-1)T}) 
= P(q^*_1 \circ \tht_{(s-1)T,(s-1)T\tht} \in \cdot |\cF_{(s-1)T}) 
=P(q^*_1 \in \cdot ).
$$
\hfill $\Box$

\vvs
By \Lem{iid}, the sequence $(|B^*_s|)_{s \in \N}$ 
is a branching process with random 
environments in the sense of Smith and 
Wilkinson \cite{SW69} (See also \cite{AtKa71a,AtKa71b}).
We write 
$$
m^*_s=\sum_{k \in \N}kq^*_s(k)
$$
for the expected number of the children in $s$ generation of the 
process $(|B^*_s|)_{s \in \N}$. 
\Lemma{N^*surv}
\bds
\item[a)]
$|B^*_s| \le |B_{sT,sT\tht}|$ for all $s \in \N$.
\item[b)]
Suppose (\ref{assP>0}). Then, 
$$
P^{\bf q} \lef( \mbox{$B^*_s \neq \epty$ for all $s \ge 1$} \ri)>0
\; \; \mbox{Q-a.s.}
$$
and hence (\ref{survP>0}) holds.
\eds
\end{lemma}
Proof: 
a) This follows easily from (\ref{branch<B}) and induction on $s$. \\
b) By \Thm{SW}(a), it is enough to show that 
$$
Q\ln m^*_1 >0\; \; \mbox{and}\; \; Q\ln {1 \over 1-q^*_1(0)} <\8.
$$
The first of the above can be seen as follows. Note that 
$B_1^*=B_{T,T\tht}$, and hence, by  (\ref{eq:defq*}),
$$
m_1^*=P^{\bf q}[|B_{T,T\tht}|]=Z_{T,T\tht}.
$$
Thus, 
$$
Q[\ln m^*_1] =Q [\ln Z_{T,T\tht}] >0.
$$
To see the second, we take $x_0,x_1,\ldots,x_T$ such that 
$x_0=0$, $x_T=\tht T$, and $|x_t-x_{t-1}|=1$, $t=1,\ldots,T$. Then,
$$
1-q^*_1(0)=P^{\bf q}(B_{T,T\tht }\neq \epty) 
\ge \prod_{t=0}^{T-1} {1-q_{t,x_t}(0) \over 2d}, 
$$
Thus,
$$
Q\ln {1 \over 1-q^*_1(0)} \le T Q\ln {1 \over 1-q_{0,0}(0)} +T\ln (2d) <\8.
$$
\hfill $\Box$

\vvs
\noindent {\bf End of the proof of \Lem{P>0}:}
The condition (\ref{KS_P>0}) reads:
$$
Q \lef[ {1 \over m^*_1}\sum_{k \ge 1}k \ln k q^*_1 (k)\ri]<\8.
$$ 
Thus, we have by \Thm{SW}(b) that
\bdnl{lnN^*}
\barray{l}
{\dps \lim_{s \ra \8}{1 \over s}\ln |B^*_s|
=
\lim_{s \ra \8}{1 \over s}\sum_{u=1}^s\ln m^*_u
=Q\lef[\ln m^*_1\ri],} \\
\mbox{a.s. on the event} \; \{B^*_s \neq \epty, \mbox{for all $s \ge 1$}\},
\earray
\edn
where the second equality comes from the law of large numbers. Since
$$
\inflim_{s \ra \8}{1 \over sT}\ln |B_{sT,sT\tht}| \ge
{1 \over T}\lim_{s \ra \8}{1 \over s}\ln |B^*_s|
$$
and  
$$
Q\lef[\ln m^*_1\ri]=
Q \ln Z_{T,T\tht}> T \ln r, 
$$
we conclude the proof of \Lem{P>0} from (\ref{lnN^*}). 
\hfill $\Box$
\subsection{Proof of \Thm{N_T}(b) and \Prop{reg=surv}}
Since $\Psi=\psi (0)$, we get (\ref{surv}) from \Lem{checkP>0}.
To show (\ref{N_T>}), it is enough to prove the following lemma.
\Lemma{701}
Suppose that
\bdnl{ass701}
P \lef( \inflim_{s \ra \8}r^{-sT}|B_{sT}| >0 \ri)>0,
\edn
for some $r>0$ and $T \in \N^*$. Then,
\bdnl{claim701}
\{\mbox{\rm survival}\} =
\{ \inflim_{t \ra \8}r^{-t}|B_t|>0\},
\; \; \; \mbox{$P$-a.s.}
\edn
\end{lemma}
Indeed, (\ref{N_T>}) follows easily from  \Lem{P>0} and \Lem{701}.
In fact, let $r=\exp ((\Psi -\e))$. Then, since $\Psi=\psi (0)$, 
\Lem{P>0} for $\tht =0$ implies (\ref{ass701}) 
and hence (\ref{claim701}) by \Lem{701}. 

\vvs
\noindent {\it Proof of \Lem{701}:} 
The following argument is adapted from \cite[page 701]{Gri83}. 
We take $u \in \{ 0,1,\ldots,T-1\}$ and fix it for a moment. 
Then, by (\ref{ass701}), the Markov property, and the shift invariance, 
we have:
$$
P \lef( \inflim_{s \ra \8}r^{-sT}|B_{sT+u}| >0 \ri)>0.
$$
This implies that
$$
P \lef( \inf_{s \ge 1}r^{-sT}|B_{sT+u}| >0 \ri)>0, 
$$
and hence that
\bds
\item[1)] \hspace{1cm}${\dps 
\del \st{\rm def}=P \lef( \inf_{s \ge 1}r^{-sT}|B_{sT+u}| >\e \ri)>0}$ 
for some $\e>0$.
\eds
We now define a series of $(\cF_{sT+u})_{s \ge 1}$-stopping times 
$0=\s_0 <\s_1 \le \s_2 \le \ldots$ as follows. 
$$
\s_1=\inf \{ s \ge 1\; ; \; 1 \le |B_{sT+u}| \le \e r^{sT} \}.
$$
Note at this point that 
\bds
\item[2)] \hspace{1cm}${\dps P(\s_1 =\8 )\ge \del }$,
\eds
thanks to 1). 
Suppose that $\s_0,\ldots,\s_\ell$ ($\ell \ge 1$) have already 
been defined. 
If $\s_\ell=\8$, we set $\s_n =\8$ for all $n \ge \ell +1$. 
Suppose that $\s_\ell<\8$. Then $B_{\s_\ell T+u} \neq \epty$ and thus, 
there is a $z \in \zd$ such that
\bds
\item[3)] $B_{\s_\ell T+u,z} \neq \epty$
\eds
and for this $z$, there is a $\lm \in \cT_{\s_\ell T+u}$ such that
\bds
\item[4)] $(z,\lm) \in B_{\s_\ell T+u}$.
\eds
Let $Z_\ell$ be the minimum, in the 
lexicographical order, of $z \in \zd$ such that 3) holds, 
and for $z=Z_\ell$, let $\lm_\ell$ be the minimum, again in the 
lexicographical order, of $\lm$ such that 4) holds. 
We now define $\s_{\ell +1}$ by:
$$
\s_{\ell +1}= \s_\ell
+\inf \{ s \ge 1 \; ;\; 1 \le |B^{\s_\ell T+u,Z_\ell,\lm_\ell}_{sT}| \le \e r^{sT} \}.
$$
(Recall that $B^{s,z,\lm}_\cdot$ denotes the $(s,z,\lm)$-branch of 
$B_\cdot$).
It is easy to see from the construction that
\bds
\item[5)] \hspace{1cm}
${\dps P(\s_\ell <\8 \; \; \mbox{i.o.})=0 }$.
\eds
Indeed, we have 
$$
P(\s_{\ell +1}<\8|\cF_{\s_\ell T+u})=P(\s_1 <\8 )\st{\scriptstyle 2)}{\le} 1-\del,
$$
and hence 
\bdnn
P(\s_{\ell +1}<\8)
&= & P(\s_\ell <\8,\; \s_{\ell +1}<\8) \\
&=&P(\s_\ell <\8, \; P(\s_{\ell +1}<\8|\cF_{\s_\ell T+u})) \\
& \le &(1-\del)P(\s_\ell <\8) \\& \le& (1-\del)^{\ell +1}
\ednn
by induction.
Then, 5) follows from the Borel-Cantelli lemma. 

By 5), we can pick a random $\ell \in \N$ such that almost surely,  
$\s_\ell <\8$ and $\s_{\ell +1} =\8$. Since $\s_\ell <\8$, we have 
$B_{t+u} \supset B^{\s_\ell T+u,Z_\ell,\lm_\ell }_{t-\s_\ell T}$ 
for all $t \ge \s_\ell T$. Note also that, on the event of survival, 
$\s_{\ell +1} =\8$ implies that
$$
|B^{\s_\ell T+u,Z_\ell,\lm_\ell}_{(s-\s_\ell )T}| 
\ge \e r^{(s-\s_\ell )T} 
\; \; \mbox{for $s \ge \s_\ell$}.
$$
Thus, on the event of survival, 
$$
|B_{sT+u}| 
\ge |B^{\s_\ell T+u,Z_\ell,\lm_\ell }_{(s-\s_\ell )T}| 
\ge \e r^{(s-\s_\ell )T}
\; \; \mbox{for $s \ge \s_\ell$}
$$
hence 
$$
\{{\rm survival}\} \st{\rm a.s.}{\sub} 
\{ \inflim_{s \ra \8}r^{-sT}|B_{sT+u}| >0\}.
$$
Since the above is true for all $u=0,1,\ldots,T-1$ and 
$$
\bigcap_{u=0}^{T-1}\{ \inflim_{s \ra \8}r^{-sT}|B_{sT+u}| >0\}
=\{\inflim_{t \ra \8}r^{-t}|B_t| >0\},
$$
we get (\ref{claim701}).
\hfill $\Box$

\vvs
\noindent {\it Proof of \Prop{reg=surv}:} 
We apply \Lem{701} to $T=1$ and $r=m$ to get the first equality. For the 
second one, 
we start to observe that the assumption implies that 
$P^{\bf q}(W_\8) >0$ with
positive $Q$-probability. By the zero-one law
(e.g., cf. \cite{CSY03}),
this event has $Q$-probability equal to 1,
and 
$$m^{-t}Z_t=P^{\bf q}[W_t] 
\to P^{\bf q}[W_\8]
\in (0,\8),
\qquad \mbox{$Q$-a.s.}
$$
Writing $m^{-t}|B_t|=m^{-t}Z_t \times Z_t^{-1}|B_t|$, we derive
the last equality.
 \hfill $\Box$
\subsection{Local survival and growth rate at all large enough times}
In this subsection, we strengthen \Lem{P>0} as follows to prepare 
the proof of \Thm{N_Tloc}(b):
\Lemma{tlarge}
Suppose $\tht \in \D \cap \Q^d$, (\ref{1-q(0)}) and that 
there exists $T \in \N^* (\tht )$ such that (\ref{assP>0}) holds. 
Then, 
\bdnl{loc_surv10}
P^{\bf q}\big( B_{t,t\tht} \neq \epty \; \; \mbox{for all large 
$t \in \N (\tht)$} \big)>0,
\; \; \mbox{$Q$-a.s.}
\edn
Suppose in addition that (\ref{KS_P>0}) holds.
Then, 
\bdnl{N_T>loc10}
P^{\bf q}\Big( |B_{t,t\tht}|  \ge e^{(\psi(\tht )-\e)t}\; \; 
\mbox{for all large $t \in \N^* (\tht)$} \Big)>0, \qquad 
\mbox{$Q$-a.s}.
\edn
\end{lemma}
Proof: 
We will prove (\ref{loc_surv10}), mentioning 
at this point that (\ref{N_T>loc10}) will be proved similarly. 
Recall notation $n(\theta)$ from (\ref{n(tht)}).
Fix $T=Kn(\theta)$ such that (\ref{assP>0}) holds.
Then, by \Lem{P>0}, 
$$
\a ({\bf q})\st{\rm def.}{=}
P^{\bf q}\big(B_{sT,sT\tht}\neq \epty \; \mbox{for all $s \ge 1$}\big) >0, 
\qquad \mbox{$Q$-a.s.}$$ 
We say that a family  $(\mu^{(i)}, i \in I)$ (with
$\mu^{(i)}\in \cT$) is {\it independent} if for all $i \neq j \in I$, 
$\mu^{(i)}$ is not an ancestor
of $\mu^{(j)}$. The reason for the terminology is that the branches
$(B^{|\mu^{(i)}|,z^{(i)},\mu^{(i)}}; i \in I)$ 
are then independent for all
choice of the $z^{(i)}$'s. 
The idea of proof is that $Q$-a.s., due to branching, 
there is a positive probability
to find a family of $K$ independent particles at times $(i+\ell K)n(\tht)$
and site $(i+\ell K)n(\tht)\theta \in \Z^d$ ($i=1, \ldots, K$). 
Then, by independence of the branches starting from these particles, 
$Q$-a.s., there is a positive probability for every such  particle
to generate a Smith-Wilkinson process in the direction $\theta$ 
which survives forever.

We now write this in details. We have, 
with $\theta_{t,x}$ the time-space shift,
\begin{eqnarray} \label{eq:pos2344}
P^{\bf q}\big(B^{t,z,\mu}_{sT,sT\tht}\neq \epty 
\; \mbox{for all $s \ge 1$}\big) 
= \alpha (\theta_{t,z} {\bf q})>0
\end{eqnarray}
We write $T(\ell,i)=(i+\ell K)n(\tht)$ to simplify the notation and let
$$E=\{ (\ell, \mu^{(1)}, \ldots, \mu^{(K)}); \ell \in \N,
\mu^{(i)} \in \cT_{T(\ell,i)}, (\mu^{(i)})_{i=1}^K \mbox{ independent}
\},$$ and, for $(\ell, \mu^{(1)}, \ldots, \mu^{(K)}) \in E$, let 
\begin{eqnarray*}
A(\ell, \mu^{(1)}, \ldots, \mu^{(K)}) &=& \bigcap_{i=1}^K
\Big\{ 
\big(T(\ell,i) \tht, \mu^{(i)}\big) \in 
B_{T(\ell,i)}
\Big\},\\
S(\ell, \mu^{(1)}, \ldots, \mu^{(K)}) &=&
\bigcap_{i=1}^K
\Big\{ 
B^{T(\ell,i),T(\ell,i)\tht, \mu^{(i)}}_{sT, sT\tht}
 \neq \epty\quad \forall s \geq 1 \Big\}.
\end{eqnarray*}
By independence,
\begin{eqnarray*}
P^{\bf q}\big( S(\ell, \mu^{(1)}, \ldots, \mu^{(K)}) \big \vert
A(\ell, \mu^{(1)}, \ldots, \mu^{(K)}) \big) 
&=& \prod_{i=1}^K \a \big( \tht_{
T(\ell,i),T(\ell,i)\tht} {\bf q} \big) \\
&>& 0, \quad \mbox{$Q$-a.s}.
\end{eqnarray*}
by (\ref{eq:pos2344}). Since 
$$ 
\bigcup_{(\ell, (\mu^{(i)})_{i=1}^K) \in E}
A(\ell, (\mu^{(i)})_{i=1}^K) \cap S(\ell, (\mu^{(i)})_{i=1}^K)
\subset
\Big\{
B_{t,t\tht} \neq \epty\; \; \mbox{for all large 
$t \in \N^* (\tht)$} \Big\},
$$
all what we need 
in order to prove our claim, is to show that
the set of environments ${\bf q}$ such that
\bdnl{P^q(cupA)>0}
P^{\bf q}\lef( \bigcup_{(\ell,  \mu^{(1)}, \ldots, \mu^{(K)} ) \in E}
A(\ell, \mu^{(1)}, \ldots, \mu^{(K)})  \ri ) >0
\edn
 has $Q$-probability 1. Observe that, since $\psi(\tht)>0$, 
we have $Q(q_{t,x}(0)+
q_{t,x}(1) \leq 1-\epsilon)>0$ for some $\epsilon \in (0,1)$, and that
$Q(q_{t,x}(0)<1)=1$ since $m>0$. 
Fix  $S=(S_i)_{i \geq 0}$ a nearest neighbor path in $\Z^d$,
with 
$S_{i n (\tht)}=i n (\tht)\tht$, 
$i=0,1,\ldots$.
With overwhelming probability as $\ell$ increases,
this path visits at least $K$ time-space
sites where branching is possible
in the first $\ell T$ steps:
$$
\lim_{\ell \to \8} Q(C_\ell) =1,$$ 
with $C_\ell \subset \Omega_{\bf q}$ given by
$$
C_\ell=\Big\{ \sum_{i=0}^{\ell T-1} {\bf 1}\{ 
q_{i,S_i}(0)+q_{i,S_i}(1) \leq 1-\epsilon\} \geq K\Big\}.
$$
Now, on $C_\ell$, the following scenario has positive 
$P^{\bf q}$-probability:
the founding ancestor starts to follow the path $S$ giving birth to (at least)
one child at each step, till the first
branching site; it splits there into (at least) 2 particles, which
continue 
to follow the path $S$ till the second
branching site; there, the first one splits into (at least) 2 particles, 
and these 3 particles continue till the next branching site, ...etc.
To write this down precisely, we introduce 
the corresponding genealogies  $(\nu^{(i)})_{i=1}^K$, i.e.,
$$
\nu^{(i)} \in {\cT}_{T(\ell,i)}
$$
composed of 1's  only except for a 2 at the time of the $i$-th branching.
Since $(\nu^{(i)})_{i=1}^K$ are independent, 
\begin{eqnarray*}
C_\ell &\subset& 
 \Big\{{\bf q} \in \Omega_q:  
P^{\bf q}[A( \ell, \nu^{(1)}, \ldots, \nu^{(K)})] \geq \epsilon^K 
(2d)^{-(\ell+1)T}
\prod_{i=0}^{(\ell+1) T-1} \big(1-q_{i,S_i}(0)\big)^K
\Big\}
\\
&\subset& 
 \Big\{{\bf q} \in \Omega_q: \; \mbox{(\ref{P^q(cupA)>0}) holds} \Big\}.
\end{eqnarray*}
Hence, 
$$Q \big( P^{\bf q}(B_{t,t\tht} \neq \epty \; \; \mbox{for all large 
$t \in \N^* (\tht)$})>0 \big) \geq Q(C_\ell),$$
which completes the proof of (\ref{loc_surv9}). The proof of
(\ref{N_T>loc9}) is totally similar, but using the second statement of 
\Lem{P>0}. 
\hfill $\Box$
\subsection{Proof of \Thm{N_Tloc}(b)}
\noindent{\it Proof of (\ref{loc_surv9}):}
By \Lem{tlarge}, we have 
$$
P^{\bf q}\big( B_{t,t\tht} \neq \epty \; \; \mbox{for all large 
$t \in \N^* (\tht)$} \big)>0,
\; \; \mbox{$Q$-a.s.}
$$
This means that $Q$-a.s. we can find a time 
$u= u({\bf q}) \in \N^* (\tht)$ and a genealogy 
$\nu= \nu ({\bf q}) \in \cT_u$ such that 
\bds \item[1)] \hspace{1cm}
${\dps 
P^{\bf q}\big( \{ (u \tht,\nu ) \in B_u\} \cap F_\nu \big)>0 \; \;}$ 
$Q$-a.s.
\eds
where
$$
F_\lm
=\{ B^{u,u\tht,\lm}_{t,(t+u)\tht} \neq \epty \; \; \mbox{for all 
$t \in \N^* (\tht)$}\}\; \; \; 
\mbox{for any $\lm \in \cT_u$.}
$$
We now take a 
nearest neighbor path $S=(S_t)_{t=0}^u$ in $\zd$ such that 
$S_t=t\tht$ for all $t \in \N (\tht )\cap [0,u]$, and fix it.
Let $\m =(1,...,1) \in \cT_u$ and 
$$
E_{\m,S} 
=\{ \mbox{The genealogy $\m$ follows $(S_t)_{t=0}^u$ and 
$K_{t,S_t,\m} \ge 1$ for all $t=0,...,u-1$}\}.
$$
Note that $E_{\m,S} $ and $F_\m$ are independent under $P^{\bf q}$, since 
$$
E_{\m,S} \in \s ((X_{t,\cdot, \m}, K_{t,\cdot, \m})_{t \le u-1}), \; \; 
F_\m \in \s ((X_{t,\cdot, \cdot}, K_{t,\cdot, \cdot})_{t \ge u}).
$$
Moreover,
\bdnn
P^{\bf q}(E_{\m,S})
&=& \prod_{t=0}^{u-1}\lef(1-q_{t,S_t} (0)\over 2d\ri)>0,\\
P^{\bf q}(F_\m)
& = & P^{\bf q}(F_\nu) 
\st{\mbox{\scriptsize 1)}}{>}0.
\ednn
Therefore, $Q$-a.s.
$$
P^{\bf q}\big( B_{t,t\tht} \neq \epty \; \; \mbox{for all 
$t \in \N (\tht)$} \big) 
\ge P^{\bf q}(E_{\m,S} \cap F_\m)=P^{\bf q}(E_{\m,S}) P^{\bf q}(F_\m)>0.
$$
 \hfill $\Box$

\vvs
\noindent{\it Proof of (\ref{N_T>loc9}):}
By \Lem{tlarge} again, we have 
$$
P^{\bf q}(|B_{t,t\tht}|  \ge e^{(\psi(\tht )-\e)t}\; \; 
\mbox{for all large $t \in \N (\tht)$})>0
\; \; \mbox{$Q$-a.s.}
$$
This, together with the similar argument as the proof of 
(\ref{loc_surv9}) implies that
$$
P^{\bf q}\lef(\inf_{t \in \N (\tht)}e^{-(\psi(\tht )-\e)t}|B_{t,t\tht}| >0\ri)>0
\; \; \mbox{$Q$-a.s.}
$$
Then, (\ref{N_T>loc9}) 
follows from the same argument as \Lem{701}. \hfill $\Box$
\subsection{Generating function of BRWRE} \label{gf}
For $q  \in \cP (\N)$, we define its {\it generating function} by
\bdnl{qhat}
{\widehat q}(s)=\sum_{k \ge 0}s^kq(k)\; \; \; s \in [0,1]
\edn 
Here and in what follows, we agree that $0^0=1$. 
For a fixed ${\bf q} \in \W_{\bf q}$ and $t \in \N$, we define 
$\Phi_t : [0,1]^{\zd} \ra [0,1]^{\zd}$ by
\bdnl{Phi_t}
\Phi_t (\xi)=(\Phi_{t,x}(\xi))_{x \in \zd}, 
\; \; 
 \Phi_{t,x}(\xi)=\sum_{y \in \zd}p(x,y){\widehat q}_{t,x}(\xi_y).
\edn
Hence $(\Phi_t)_{t \in \N}$ is a sequence of 
i.i.d. random maps on the probability 
space $(\W_{\bf q}, \cF_{\bf q}, Q)$.
\Lemma{xi^B}
For $\xi \in [0,1]^{\zd}$ and $t \in \N^*$, 
\bdnl{xi^B|cF_t}
P^{\bf q}[\xi^{B_t}| \cF_{t-1}] =\Phi_{t-1}(\xi)^{B_{t-1}},
\edn
with the notation
$$
\xi^{B_t} =\prod_{x \in \zd}\xi_x^{|B_{t,x}|} \in [0,\8).
$$
As a consequence, 
\bdnl{xi^B}
P^{\bf q}[\xi^{B_t}] 
=\Phi_{0,0} \circ \Phi_1 \circ \ldots \circ \Phi_{t-1}(\xi).
\edn
\end{lemma}
Proof: It is enough to show (\ref{xi^B|cF_t}). We begin by writing:
\bdnn
\xi^{B_t}=\prod_{y \in \zd}\xi_y^{|B_{t,y}|}
&=&\prod_{y \in \zd}\prod_{(x,\n) \in B_{t-1}}
\xi_y^{K_{t-1,x,\n}\del_y (X_{t-1,x,\n})}
=\prod_{(x,\n) \in B_{t-1}}
\prod_{y \in \zd}\xi_y^{K_{t-1,x,\n}\del_y (X_{t-1,x,\n})} \\
&=&\prod_{(x,\n) \in B_{t-1}}
\prod_{y \in \zd}\lef( 1+(\xi_y^{K_{t-1,x,\n}}-1)\del_y (X_{t-1,x,\n}) \ri) \\
&=&\prod_{(x,\n) \in B_{t-1}}
\lef(1+\sum_{y \in \zd}(\xi_y^{K_{t-1,x,\n}}-1)\del_y (X_{t-1,x,\n}) \ri) \\
&=&\prod_{(x,\n) \in B_{t-1}}
\sum_{y \in \zd}\xi_y^{K_{t-1,x,\n}}\del_y (X_{t-1,x,\n}) 
\ednn
where on the third line, we have used that
\bds \item[1)] \hspace{1cm}
${\dps 
\prod_{y \in Y}(1+x_y)=1+\sum_{A \sub Y}\prod_{y \in A}x_y}$ \quad
for any finite 
set $Y$ and $(x_y)_{y \in Y} \in \R^Y$
\eds 
(the terms with $|A| \ge 2$ on the right-hand-side of 1) 
vanishes in our application). Since 
$$
P^{\bf q} \lef[ 
\sum_{y \in \zd}\xi_y^{K_{t-1,x,\n}}\del_y (X_{t-1,x,\n}) | \cF_{t-1} \ri]
=\Phi_{x,t-1}(\xi),
$$
We get 
$$
P^{\bf q}[\xi^{B_t}| \cF_{t-1}]
=\prod_{(x,\n) \in B_{t-1}}\Phi_{x,t-1}(\xi)
=\Phi_{t-1}(\xi)^{B_{t-1}}.
$$
\hfill $\Box$

\vvs
Proof of \Prop{surv01}: 
We note that the map $\Phi_t$ ($t \in \N$) 
has the following continuity property:
\bdnl{Phi_tC}
\xi,\xi_n \in [0,1]^{\zd}\; \; \lim_n\xi_n=\xi
\; \; \Longrightarrow \; \; \lim_n\Phi_t(\xi_n)=\Phi_t(\xi),
\edn
where the limits are coordinatewise. We have 
also that
\bdnl{Phi(1)}
\Phi_{0,0} \circ \Phi_1 \circ \ldots \circ \Phi_{t-1}(\xi)=1 
\; \; \Longleftrightarrow \; \; 
\xi_x=1\; \; \mbox{for all $x \in \zd$ with $P_S (S_t=x)>0$}.
\edn
Let $B^x_{\cdot}$ be the BRWRE 
starting from one particle from time-space 
$(0,x) \in \N \times \zd$. Then, for the zero-field 
$\xi$ (i.e., $\xi_\cdot \equiv 0$), 
\bds \item[1)] \hspace{1cm}
$\del_{t,x} ({\bf q})
\st{\rm def}{=} P^{\bf q}(B^x_t =\epty)=P^{\bf q}[\xi^{B^x_t}]
\st{\mbox{\scriptsize (\ref{xi^B})}}{=}
\Phi_{0,x} \circ \Phi_1 \circ \ldots \circ \Phi_{t-1}(\xi).$
\eds
Since $\del_{t,x} ({\bf q})$ is decreasing in $t$, 
the limit $\del_x ({\bf q})=\lim_{t \ra \8}\del_{t,x} ({\bf q})$ 
exists for all 
$x \in \zd$ and we define the random field of extinction probabilities:
\bdnl{del(q)}
\del({\bf q})=(\del_x ({\bf q}))_{x \in \zd}.
\edn
Note that $\del_0 ({\bf q})$ is the extinction 
probability of $B_t$:
\bdnl{del_0(q)}
\del_0 ({\bf q})=1-P^{\bf q}({\rm survival}).
\edn
Now, we have for any $u \in \N^*$ that
\bdmn
\del ({\bf q})
& \st{\mbox{\scriptsize 1)}}{=} &
\lim_{t \ra \8}\Phi_0 \circ \Phi_1 \circ \ldots \circ \Phi_{t-1}(\xi)\nn \\
&\st{\mbox{\scriptsize (\ref{Phi_tC})}}{=}&
\Phi_0 \circ \Phi_1 \circ \ldots \circ \Phi_{u-1}
\lef( \lim_{t \ra \8}\Phi_u \circ \ldots \circ \Phi_{u+t-1}\ri) \nn \\
& \st{\mbox{\scriptsize 1)}}{=} &
 \Phi_0 \circ \Phi_1 \circ \ldots \circ \Phi_{u-1}
(\del (\tht_{u,0}{\bf q})), \label{del=}
\edmn
where $\tht_{u,y}$ denotes the shift: 
$\tht_{u,y}{\bf q}=(q_{\cdot+u,\cdot+y})$. This and 
(\ref{Phi(1)}) imply that 
$$
\{ {\bf q} \; ; \; \del_0 ({\bf q})=1\}=\bigcap_{u \ge 1}
\{ {\bf q} \; ; \; \del_x (\tht_u{\bf q})=1\; 
\mbox{for all $x \in \zd$ with $P_S (S_t=x)>0$}\}.
$$
Since the right-hand-side is a tail event, it is 
$Q$-trivial by Kolmogorov's zero-one law. \hfill $\Box$

\vvs
In analogy with GW and SW processes, 
the extinction probability of BRWRE can be characterized in 
terms of the functional equation involving the generating function. 
To explain it, we introduce the coordinatewise order in $[0,1]^{\zd}$.
Let $\xi,\h \in [0,1]^{\zd}$. We write $\xi \le \h$ 
when $\xi_x \le \h_x$ for all $x \in \zd$. 
A function $F: [0,1]^{\zd} \ra [0,1]^{\zd}$ is called 
{\it increasing} if it preserves this order. 
\Corollary{minisol}
The random field $\del({\bf q})$ defined by 
(\ref{del(q)}) 
is the minimal among all 
random fields $\xi ({\bf q})$ such that 
\bdnl{fixed}
\xi ({\bf q})=\Phi_0 (\xi (\tht_{1,0}{\bf q})).
\edn
Moreover, the extinction probability 
$\del_0 ({\bf q})$ (cf. (\ref{del_0(q)})) is the 
minimal $[0,1]$-valued function 
$\xi_0 ({\bf q})$ of ${\bf q} \in \cP(\N)^{\N \times \zd}$ such that
\bdnl{fixed_0}
\xi_0 ({\bf q})=\Phi_{0,0} \lef( (\xi (\tht_{1,x}{\bf q})_{x \in \zd}) \ri).
\edn
\end{corollary}
Proof:
By setting 
$u=1$ in (\ref{del=}), we see that $\del({\bf q})$ is indeed 
a solution of (\ref{fixed}). On the other hand, any solution 
$\xi ({\bf q})$ of (\ref{fixed}) satisfies:
\bdnn
\xi ({\bf q})
&=&\Phi_0 (\xi (\tht_{1,0}{\bf q}))
=\Phi_0 \circ \Phi_1 (\xi (\tht_{2,0}{\bf q}))=\ldots \\
&=& \lim_{t \ra \8}\Phi_0 \circ \ldots 
\circ \Phi_{t-1} (\xi (\tht_{t,0}{\bf q})).
\ednn
Comparing this with the definition of $\del ({\bf q})$, we 
see that $\del ({\bf q}) \le \xi ({\bf q})$, since $\Phi_u$, $u \in \N$ 
are increasing. These prove the first half of the corollary.

A function $\xi_0: \cP(\N)^{\N \times \zd} \ra [0,1]$ solves 
(\ref{fixed_0}) if and only if  the random field defined by 
$\xi ({\bf q})=\lef( \xi_0 (\tht_{0,x}{\bf q})\ri)_{x \in \zd}$ 
solves (\ref{fixed}). 
Thus, the second half of the corollary comes down to the first half.
\hfill $\Box$

\vvs
We next turn to the proof of (\ref{SW<surv<GW}).
\Lemma{GWgen}
Let $q(k)=Q[q_{0,0}(k)]$, $k \in \N$. 
Then, for $s \in [0,1]$, 
\bdnl{GWgen}
P[s^{|B_t|}] \ge 
(\underbrace{{\widehat q} \circ {\widehat q} 
\circ \ldots  \circ {\widehat q}}_t)(s).
\edn
\end{lemma}
Proof:
It follows from 
(\ref{xi^B|cF_t}) and H\"older's inequality that
$$
P[s^{|B_t|}|\cF_{t-1}]
=\prod_{x \in \zd}Q[{\widehat q}_{0,0}(s)^n]_{n=|B_{t-1,x}|}
\ge \prod_{x \in \zd}Q[{\widehat q}_{0,0}(s)]^{|B_{t-1,x}|}
={\widehat q} (s)^{|B_{t-1}|}.
$$
Thus, (\ref{GWgen}) follows by iteration.
\hfill $\Box$
\Lemma{SWgen}
For $s \in [0,1]$, 
\bdnl{SWgen}
P[s^{|B_t|}] \le Q[
{\widehat q}_{0,0} \circ {\widehat q}_{1,0} \circ \ldots  \circ 
{\widehat q}_{t-1,0}(s)].
\edn
\end{lemma}
Proof:
We will prove by induction that for $1 \le u \le t$, 
\bds \item[1)] \hspace{1cm} ${\dps 
P[s^{|B_t|}|\cF_{t-u}] \le Q[
\lef({\widehat q}_{0,0} \circ {\widehat q}_{1,0} \circ \ldots  \circ 
{\widehat q}_{u-1,0}\ri)(s)^n]_{n=|B_{t-u}|}.
}$ \eds
This in particular proves (\ref{SWgen}) ($u=t$). It follows from 
(\ref{xi^B|cF_t}) and H\"older's inequality that
$$
P[s^{|B_t|}|\cF_{t-1}]
=\prod_{x \in \zd}Q[{\widehat q}_{0,0}(s)^n]_{n=|B_{t-1,x}|}
\le Q[{\widehat q}_{0,0}(s)^n]_{n=|B_{t-1}|},
$$
which proves 1) for $u=1$. Now assume that 1) with $u$ replaced by 
$u-1$ ($u \ge 2$) is true and let $A \in \cF_{t-u}$. Then, by the 
assumption and Fubini's theorem,
\bdnn
E[s^{|B_t|}{\bf 1}_A] 
& \le &
P\lef[ {\bf 1}_AQ[({\widehat q}_{0,0} \circ {\widehat q}_{1,0}
\circ  \ldots \circ {\widehat q}_{u-2,0})(s)^n]_{n=|B_{t-u+1}|} \rig] \\
& = &
Q\lef[ P[ {\bf 1}_A ({\widehat q}_{t-u+1,0} \circ {\widehat q}_{t-u+2,0} 
\circ \ldots  \circ {\widehat q}_{t-1,0})(s)^{|B_{t-u+1}|} ] \ri].
\ednn
By 1) with $u=1$, we have that 
\bdnn
\lefteqn{P[ {\bf 1}_A ({\widehat q}_{t-u+1,0} \circ {\widehat q}_{t-u+2,0}
\circ \ldots  \circ {\widehat q}_{t-1,0})(s)^{|B_{t-u+1}|} ]} \\
&\le &P\lef[ {\bf 1}_A\int Q (dq_{t-u,0})
({\widehat q}_{t-u,0} \circ {\widehat q}_{t-u+1,0} 
\circ  \ldots \circ {\widehat q}_{t-1,0})(s)^{|B_{t-u}|} \rig].
\ednn
Here, $\int Q (dq_{t-u,0})$ means that we only integrate 
$q_{t-u, 0}$, with all the other $q_{\cdot, \cdot}$ fixed.
Combining these, we arrive at 
$$
E[s^{|B_t|}{\bf 1}_A] 
\le P\lef[ {\bf 1}_AQ[ 
({\widehat q}_{t-u,0} \circ {\widehat q}_{t-u+1,0}
\circ  \ldots \circ {\widehat q}_{t-1,0})(s)^n]_{n=|B_{t-u}|} \rig],
$$
which proves 1). \hfill $\Box$

\vvs
\noindent {\it Proof of (\ref{SW<surv<GW}):}
Let $(B^{\rm GW}_t)_{t \in \N}$ and $(B^{\rm SW}_t)_{t \in \N}$ 
be the Galton-Watson and Smith-Wilkinson processes we are interested in.
Then, for $s \in [0,1]$, 
\bdnn
P[s^{|B^{\rm GW}_t|}] 
&=&(\underbrace{{\widehat q} \circ {\widehat q} 
\circ \ldots  \circ {\widehat q}}_t)(s),\\
P[s^{|B^{\rm SW}_t|}] 
&=&
Q[{\widehat q}_{0,0} \circ {\widehat q}_{1,0} \circ \ldots  \circ 
{\widehat q}_{t-1,0}(s)].
\ednn
The former is well-known and for the latter, see \cite[Theorem 2.1]{SW69}.
For $s=0$, we have 
$$
P({\rm survival})=1-\lim_{t \ra \8}P[s^{|B_t|}].
$$
and similar formulae for $\s^{\rm GW}$ and $\s^{\rm SW}$. 
Therefore, (\ref{SW<surv<GW}) follows from \Lem{GWgen} and \Lem{SWgen}. 

\hfill $\Box$

\SSC{Appendix}
\subsection{Proof of \Prop{press}(a)} \label{ppressb}
If $m_{0,0} \equiv m$ a.s., then 
$\psi (\tht)=\ln m -I_S (\tht)$. Thus, we assume 
 that $m_{0,0} \not\equiv m$ a.s.

We start by proving the following. 
\bds
\item[a)]
For $\tht \in \D \cap \Q^d$ and $t \in \N^* (\tht)$, 
\bdnl{Ztht>}
\psi (\tht) \ge \frac{1}{t} Q[\ln Z_{t,t\tht}]
\ge Q[\ln m_{0,0}]+\frac{1}{t} \ln P_S( S_t=t\tht ).
\edn
Moreover, 
the second inequality is an equality if and only if $t =1$ or  
$\tht \in \{ \pm e_i\}_{i=1}^d$.
\item[b)]
For $t \ge 1$, 
\bdnl{Z>}
\psi (0) \ge  \frac{1}{t} Q[\ln Z_t] \ge Q[\ln m_{0,0}].
\edn 
Moreover, 
the second inequality is an equality if and only if $t=1$.
\eds
These can be seen as follows:
\bdnn
Z_{t,t\tht } 
&=& P_S[\prod^{t-1}_{u=0}m_{u,S_u}| S_t=t\tht ] 
P_S( S_t=t\tht ) \\
&\st{\mbox{\scriptsize Jensen}}{\ge}&
\exp \lef( 
\sum^{t-1}_{u=0}P_S[\ln m_{u,S_u}|S_t=t\tht] \ri)P_S( S_t=t\tht ).
\ednn
Note that, if $t \neq 1$ and $\tht \not\in \{ \pm e_i\}_{i=1}^d$ 
the random variable 
$S_1$ under $P_S (\; \cdot \; | S_t=t\tht)$ is not a constant.
Thus, the Jensen inequality above is strict on the event:
$$
\{ {\bf q} \; ; \; m_{1,e} \neq m_{1,e'} \; 
\mbox{for some $e, e' \in \zd$ with $|e|=|e'|=1$}\},
$$
which has positive $Q$-probability by the assumption. 
By taking logarithm, and then $Q$-expectation, we get the 
second inequality of (\ref{Ztht>}). 
On the other hand, it is not difficult to prove that 
the sequence $\{ Q[\ln Z_{t,t\tht}]\}_{t \in \N^*}$ is 
superadditive (e.g.,\cite[page 720, Proof of Proposition 2.5]{CSY03}) 
and hence 
$$
\psi (\tht)=\sup_{t \ge 1}\frac{1}{t}Q[\ln Z_{t,t\tht}].
$$ 
This completes the proof of of (\ref{Ztht>}).
The proof of (\ref{Z>}) is similar.

By letting $t \ra \8$ in (\ref{Ztht>}), we get 
$$
\psi (\tht)\ge Q[\ln m_{0,0}]-I_S(\tht ).
$$
Finally, we prove 
$$
\psi (\tht) \le \ln m-I_S (\tht).
$$
We have 
\bdnn
Q [\ln Z_{t,t\tht } ] &=& 
Q [\ln P_S[\prod^{t-1}_{u=0}m_{u,S_u}: S_t=t\tht ]] \\
&\st{\mbox{\scriptsize Jensen}}{\le}&
\ln Q [ P_S[\prod^{t-1}_{u=0}m_{u,S_u}: S_t=t\tht ]] 
= t \ln m+\ln P_S(S_t=t\tht )
\ednn 
We get the desired bound by 
dividing the above inequality by $t$ and then letting 
$t \ra \8$.
\hfill $\Box$
\subsection{Proof of \Prop{press}(b)} \label{ppressc}
Let $(\W, \cG, Q)$ be a probability space and 
$$
\{ \epty, \W \}=\cG_0 \sub \cG_1 \sub \ldots 
$$
be an increasing series of sub $\s$-fields of $\cG$. For 
$X \in L^1 (Q)$, we write $Q^{\cG_{j}}[X]=Q[X| \cG_j]$.
We use a concentration inequality in the following Lemma. 
\Lemma{lem-concentration}
 Suppose that $X \in L^1(Q)$ is $\cG_n$-measurable for some $n$ 
and that there exist $\del \in (0,\8)$, $A \in [0,\8)$, 
$X_1,\ldots 
,X_n \in L^1(Q)$ such that 
\bdnl{X-Q[X]1}
Q^{\cG_{j-1}}[X_j]=Q^{\cG_{j}}[X_j],\; \; \mbox{and}\; \; 
Q^{\cG_{j-1}}[\exp(\del |X-X_j|)] \le A
\edn
for all $j=1,\ldots ,n$. 
Then, with $B=2\sqrt{6}A^2/\del^2$, 
\bdnl{X-Q[X]2}
Q (|X-Q[X]| \ge \e n)  
 \le  2\exp ( -B\e^2n/4)\; \; \; 
\mbox{for all  $\e \in [0, B\del ]$}. 
\edn
\end{lemma}
Proof: 
We consider a  sequence $D_j=Q^{\cG_{j}}[X]-Q^{\cG_{j-1}}[X]$.  We first 
observe that 
\bds
\item[1)] \hspace{1cm}
$ Q^{\cG_{j-1}}[e^{\del |D_j|}] \le A^2$ 
for all $j=1,\ldots ,n$. 
\eds
Since $Q^{\cG_{j-1}}[X_j]=Q^{\cG_{j}}[X_j]$, we have 
\bdnn
|D_j| 
& \le & |Q^{\cG_{j}}[X-X_j]|+|Q^{\cG_{j-1}}[X-X_j]| \\
& \le & Q^{\cG_{j}}[Y_j]+Q^{\cG_{j-1}}[Y_j], \; \; 
\mbox{with $Y_j=|X-X_j|$.}
\ednn
It follows from Jensen inequality that
$$
e^{\del Q^{\cG_{j-1}}[Y_j]} \le Q^{\cG_{j-1}}[e^{\del Y_j}] \le A.
$$
Similarly, 
$$
Q^{\cG_{j-1}}[e^{\del Q^{\cG_{j}}[Y_j]}] 
\le Q^{\cG_{j-1}}[Q^{\cG_j}[e^{\del Y_j}] ] 
= Q^{\cG_{j-1}}[e^{\del Y_j} ] \le A.
$$
These imply 1) as follows:
$$
Q^{\cG_{j-1}}[e^{\del |D_j|}] \le e^{\del Q^{\cG_{j-1}}[Y_j]}
Q^{\cG_{j-1}}[e^{\del Q^{\cG_{j}}[Y_j]}] \le A^2.
$$
We now infer from 1) that
\bds
\item[2)] \hspace{1cm}
${\dps Q^{\cG_{j-1}}[e^{\a D_j}] 
\le e^{B\a^2}}$ 
for all $\a \in [-\del /2, \del /2]$ and $j=1,\ldots ,n$. 
\eds
Observe that
$$
\mbox{$\frac{1}{4!}$}Q^{\cG_{j-1}}[|D_j|^4] = 
\mbox{$\frac{1}{4!}$}Q^{\cG_{j-1}}[\del^4|D_j|^4]/\del^4 \le 
Q^{\cG_{j-1}}[e^{\del|D_j|}]/\del^4 \le A^2/\del^4
$$
and hence that
$$
Q^{\cG_{j-1}}[|D_j|^2e^{\del  |D_j|/2}] \le 
Q^{\cG_{j-1}}[|D_j|^4]^{1/2}
Q^{\cG_{j-1}}[e^{ \del | D_j|}]^{1/2}
\le 2\sqrt{6}A^2/\del^2=B.
$$
Since 
$e^{x} \le 1 +x +|x|^2e^{|x|}/2$ for all $x \in \R$, we get
$$
Q^{\cG_{j-1}}[e^{\a D_j}] 
\le 1 +\a^2B/2 \le \exp (B\a^2)
$$
Finally, since $X-Q[X]=D_n+\ldots +D_1$, it follows from 2) that 
\bds
\item[3)] \hspace{1cm}
$Q [\exp( \a (X-Q[X]))]  \le \exp ( B\a^2n)\; \; \; 
\mbox{for all  $\a \in [-\del /2,\del /2]$}$
\eds
via a simple iterative procedure. 
To see (\ref{X-Q[X]2}), 
we take $\a=\frac{\e}{2B} \le \frac{\del}{2}$. Then, by 
Chebychev's inequality and 3), 
\bdnn
Q (|X-Q[X]| \ge \e n) 
& = & Q (\a |X-Q[X]| \ge \a \e n) \\
& \le & 2\exp ( (B\a^2 -\a \e )n)=2\exp ( -B\e^2n/4).
\ednn
\hfill $\Box$

\vvs 
\noindent {\bf Proof of \Prop{press}(b):}
We fix an arbitrary integer $t$ and set 
$\cG_j =\s (m_{k,\cdot}: k \le j)$ ($j=1,...,t$).
To prove (\ref{exp_con_loc}) and (\ref{exp_con}) at a stroke, 
we define
$$
X=\ln P_S [\z_t f(S_t)],
$$
where $f: \zd \ra [0, \8)$ is such that $P_S [f(S_t)] > 0$. 
$X=\ln Z_{t,t\tht}$ for $f (x)=\del_{x, t\tht}$ and, 
$X=\ln Z_t$ for $f (x)\equiv 1$. 
We will prove that $X$ satisfies (\ref{X-Q[X]2}) 
with some $B$, independent of the choice of $f$.
To do so, let
$$
X_j=\ln P_S [\z_{t,j} f(S_t)],\; \; j=1,...,t,
$$
where
$$
\z_{t,j}=\prod_{0 \le k \le t-1 \atop k \neq j}m_{k,S_k}.
$$
We then check (\ref{X-Q[X]1}).
Since $X_j$ does not depend on $m_{\cdot, j}$, we have 
$Q^{\cG_{j-1}}[X_j]=Q^{\cG_{j}}[X_j]$. Note on the other hand that
for $\del \in \R \bsh (0,1)$, 
\bdnn
\exp (\del(X-X_j))
 =  \lef( {P_S [\z_t f(S_t)] \over P_S [\z_{t,j} f(S_t)]} \ri)^\del
& = &  \lef( {P_S [\z_{t,j}m_{j,S_j} f(S_t)] 
\over P_S [\z_{t,j} f(S_t)]} \ri)^\del\\
& \st{\mbox{\scriptsize Jensen}}{\le} & 
{P_S [\z_{t,j}m_{j,S_j}^\del f(S_t)] \over P_S [\z_{t,j} f(S_t)]}
\ednn
Now, by taking the conditional expectation given 
$\cG_j^\pri \st{\rm def.}{=}\s [m_{k,\cdot}\; ; \; k \neq j]$, 
$$
Q^{\cG_j^\pri}[\exp (\del (X-X_j))]  \le  Q[m_{0,0}^\del].
$$
This, together with $\cG_{j-1} \sub \cG_j^\pri$, implies
$$
Q^{\cG_{j-1}}[\exp (\del(X-X_j))] \le Q[m_{0,0}^\del].
$$
By applying this for $\del =\pm 1$, we get 
$$
Q^{\cG_{j-1}}[\exp (|X-X_j|)] \le A:=Q[m_{0,0}+m_{0,0}^{-1}].
$$
\hfill $\Box$

\bigskip
{\bf Acknowledgements:} 
We warmly thank
the anonymous referee for his careful reading and numerous 
suggestions.
\small

Francis Comets \\
Universit\'e Paris Diderot -- Paris 7\\
Math\'ematiques, Case 7012, B\^atiment Chevaleret\\
75205 PARIS CEDEX 13 \\
France \\
{\tt comets@math.jussieu.fr }
\vvs

Nobuo Yoshida \\
Division of Mathematics, Graduate School of Science, \\
Kyoto University,\\
Kyoto 606-8502, \\
Japan, \\
{\tt nobuo@math.kyoto-u.ac.jp}

\begin{thebibliography}{99}
\bibitem{AtKa71a}
Athreya, K.; Karlin, S.:  
Branching processes with random environments. I. Extinction probabilities. 
{\it Ann. Math. Statist.} {\bf  42}  (1971), 1499--1520. 
\bibitem{AtKa71b}
Athreya, K.; Karlin, S.:  
Branching processes with random environments. II. Limit theorems.
{\it Ann. Math. Statist.} {\bf  42}  (1971), 1843--1858.  
\bibitem{AtNe72}
Athreya, K. and Ney, P. (1972): {\it Branching Processes},
Springer Verlag New York.
\bibitem{BG90}
Bezuidenhout, C.and  Grimmett, G.: 
The critical contact process dies out, 
Ann. Prob. {\bf 18}, No.2, (1990), 1462--1482.
\bibitem{Big79}
Biggins, J.:
Growth rates in the branching random walk.
{\it Z. Wahrsch. Verw. Gebiete} {\bf 48} (1979) 17--34. 
\bibitem{Bir03}
Birkner, M.(2003): Particle systems with locally dependent
branching: long-time behaviour, genealogy and critical parameters.
PhD thesis, Johann Wolfgang Goethe-Universit\"at, Frankfurt.
\bibitem{BGK05}
 Birkner, M., Geiger, J., Kersting G.(2005):
Branching processes in random environment-- a view on critical and
subcritical cases. Interacting stochastic systems,  269--291,
Springer, Berlin.
\bibitem{CaHu03}
Carmona, P., Hu Y.: 
Fluctuation exponents and large deviations 
for directed polymer in a random environment,
{\it Stoch. Proc. Appl.} {\bf 112}, (2004), 285-308. 
\bibitem{CP07}
Comets, F., Popov, S.: 
Shape and local growth for multidimensional branching random walks in random environment.  
{\it ALEA Lat. Am. J. Probab. Math. Stat.} {\bf  3}  (2007), 273--299
\bibitem{CSY03}
Comets, F., Shiga, T., Yoshida, N. Directed Polymers in
Random Environment: Path Localization and  Strong Disorder, {\it
Bernoulli} {\bf 9} (2003) 705--723.
\bibitem{CSY04}
Comets, F., Shiga, T., Yoshida, N. (2004) Probabilistic analysis
of directed polymers in random environment: a review, Advanced
Studies in Pure Mathematics, {\bf 39}, 115--142.
\bibitem{CV04}
Comets, F., Vargas, V. 
Majorizing multiplicative cascades for directed polymers
in random media. {\it ALEA Lat. Am. J. Probab. Math. Stat.}  {\bf 2}
(2004)
267--277.
\bibitem{Dur95}
Durrett, R. (2005):``Probability--Theory and Examples", 3rd Ed.,
Brooks/Cole--Thomson Learning.
\bibitem{GMPV09}
Gantert, N., M\"uller, S.,  Popov, S., Vachkovskaia, M.
 Survival of Branching Random Walks in Random Environment.
{\it  J. Theoret. Probab.}, to appear
\bibitem{GdH91}
Greven, A., den Hollander, F.:
Branching random walk in random environment: phase transitions for local 
and global growth rates.
{\it Probab. Theory Related Fields} {\bf 91} (1992) 195--249. 
\bibitem{Gri83}
Griffeath, D.:
The Binary Contact Path Process, 
{\it Ann. Probab.} { 11} (1983) 692-705.
\bibitem{GH02}
Grimmett, Geoffrey; Hiemer, Philipp:
Directed percolation and random walk. 
In and out of equilibrium (Mambucaba, 2000), 273--297, 
Progr. Probab., 51, Birkha\"user,
 (2002).
\bibitem{HuYo07}
Hu, Y., Yoshida, N. :
Localization for Branching Random Walks in Random Environment,
{\it Stoch. Proc. Appl.}
{\bf 119} (2009) 1632--1651.
\bibitem{KeSti66a}
Kesten, H.; Stigum, B. P. :
A limit theorem for multidimensional Galton-Watson processes.  
{\it Ann. Math. Statist.} {\bf  37}  (1966) 1211--1223. 
\bibitem{La09} 
Lacoin, H. :
New bounds for the free energy of directed polymers 
in dimensions 1+1 and 1+2, preprint.
\bibitem{LW08} 
Liu, Q.; Watbled, F.:
Exponential inequalities for martingales and asymptotic 
properties of the free energy of directed polymers in a random environment,
preprint, 2008, to appear in Stoch. Proc. Appl. 
\bibitem{Mu08}
M\"uller, S.:  A criterion for transience of multidimensional branching random walk in random environment.  
{\it Electron. J. Probab.}  {\bf 13}  (2008) 1189--1202.
\bibitem{Na09} Nakashima, M:
Almost Sure Central Limit Theorem for 
Branching Random Walk in Random Environment, 
preprint 2008.
\bibitem{Rev94}
R\'ev\'esz, P. (1994): ``Random Walks of Infinitely Many
Particles" World Scientific.
\bibitem{Roc70}
Rockafeller, R. T. (1970):
Convex Analysis, Princeton University Press.
\bibitem{Sm68}
Smith, W.:
Necessary conditions for almost sure extinction of a branching process 
with random environment.
{\it Ann. Math. Statist.} {\bf 39}  (1968) 2136--2140. 
\bibitem{SW69}
Smith, W., L.; Wilkinson, W., E.: On branching processes in
random environments. {\it Ann. Math. Statist.} {\bf  40} (1969) 814--827.
\bibitem{SoZh96}
Song, R. and Zhou, X. Y. : 
A remark on diffusion on directed polymers in random environment, 
{\it J. Stat. Phys.} {\bf 85} 277--289 (1996).
\bibitem{Yo08a} 
Yoshida, N.: Central Limit Theorem for 
Branching Random Walk in Random Environment, 
{\it Ann. Appl. Proba.} {\bf 18} (2008) 1619--1635.
\bibitem{Yo08b} 
Yoshida, N.:
Phase Transitions for the Growth Rate of Linear Stochastic Evolutions, 
{\it J. Stat. Phys.} {\bf 133} 1033--1058 (2008). 
\end{thebibliography}
\end{document}